\documentclass[a4paper,11pt]{article}
\usepackage[utf8]{inputenc}
\usepackage[T1]{fontenc}
\usepackage[english]{babel}
\usepackage{xcolor,graphicx,tikz}
\usepackage{amsthm,amsmath,amssymb}
\usepackage[top=2.5 cm,bottom=2 cm,left=2.5 cm,right=2 cm]{geometry}
\usepackage[utf8]{inputenc}
\usepackage{mathrsfs}
\usepackage{bbm}
\usepackage{esint}
\usepackage{graphicx}
\usepackage{hyperref}\hypersetup{colorlinks=true, citecolor=blue}

\usepackage{color}
\usepackage{amssymb}
\usepackage{amsmath,amscd}
\usepackage{mathtools}
\mathtoolsset{showonlyrefs,showmanualtags}

\newcounter{stepnb}

\theoremstyle{definition}

\theoremstyle{definition}
\newtheorem{question}{Question}

\newcommand{\RR}{\mathbb{R}}

\newcommand{\ee}{\varepsilon}

\graphicspath{
{../Figures/}
{../Figures/Figures_New_Schemes/}
}
\usepackage{authblk}
\title{On the role of numerical viscosity in the study \\ of the local limit of nonlocal conservation laws  }
\author[1]{Maria Colombo}
\author[2]{Gianluca Crippa}
\author[3]{Marie Graff}
\author[4]{Laura V. Spinolo}
\affil[1]{ EPFL SB, Station 8, CH-1015 Lausanne, Switzerland.
E-mail: maria.colombo@epfl.ch}
\affil[2]{Departement Mathematik und Informatik,
Universit\"at Basel, Spiegelgasse 1, CH-4051 Basel, Switzerland. Email: gianluca.crippa@unibas.ch}
\affil[3]{Department of Mathematics, University of Auckland, New Zealand. Email: marie.graff@auckland.ac.nz}
\affil[4]{IMATI-CNR, via Ferrata 5, I-27100 Pavia, Italy. Email: spinolo@imati.cnr.it}
\date{}                 

\begin{document}
\maketitle

\begin{abstract}
         We deal with the numerical investigation of the local limit of nonlocal conservation laws. Previous numerical experiments suggest convergence in the local limit. However, recent analytic results state that (i) in general convergence does not hold because one can exhibit counterexamples; (ii)~convergence can be recovered provided viscosity is added to both the local and the nonlocal equations.  Motivated by these analytic results, we investigate the role of numerical viscosity in the numerical study of the local limit of nonlocal conservation laws. In particular, we show that the numerical viscosity of Lax-Friedrichs type schemes jeopardizes the reliability of the numerical scheme and erroneously detects convergence in cases where convergence is ruled out by analytic results. We also test Godunov type schemes, less affected by numerical viscosity, and show that in some cases they provide more reliable results. 
\end{abstract}

\tableofcontents

%\newpage
\section{Introduction}
\label{s:Intro}
\subsection{Theoretical framework}
We consider nonlocal conservation laws in the form
\begin{equation}
\label{e:nonlocal}
          \partial_t \rho + \partial_x \big[ \rho \  b(\rho \ast \eta) \big] = 0, 
\end{equation} 
where the unknown is $\rho: [0, + \infty) \times \RR  \to \RR$, $b: \RR \to \RR$ is a given Lipschitz continuous function and $\eta: \RR \to \RR$ is a smooth convolution kernel satisfying  
\begin{equation}
\label{e:eta}
      \eta \in C_c^\infty (\RR), \quad \eta (x) = 0 \; \text{if $|x| \ge 1$}  , \quad \eta \ge 0, \quad \int_{\RR} \eta (x) dx =1. 
\end{equation} 
In recent years, nonlocal conservation laws have been used to model, among others, sedimentation~\cite{Sedimentation}, pedestrian~\cite{ColomboGaravelloMercier} and vehicular~\cite{BlandinGoatin,ChiarelloGoatin} traffic. In particular, in the case of traffic models 
$\rho$ represents the density of agents (cars, pedestrians) and $b$ their speed. The convolution term models the fact that 
drivers and pedestrians decide their velocity based on the density of agents around them. Loosely speaking, the radius of the support of $\eta$ represents the visual range of drivers and pedestrians. Existence and uniqueness results for the Cauchy problem obtained by coupling~\eqref{e:nonlocal} with an initial datum have been obtained in several works, see for instance~\cite{BlandinGoatin,ColomboGaravelloMercier,CLM}. 

In this work we deal with the numerical investigation of the local limit. More precisely, we consider a parameter $\ee >0$ and we rescale $\eta$ by setting $\eta_\ee (x) : =\ee^{-1} \eta \left( x /\ee \right)$, in such a way that, when $\ee \to 0^+$,
 $\eta_\ee$ converges weakly-$^\ast$ in the sense of measures to the Dirac delta. We fix an initial datum $\bar  \rho: \RR \to \RR$ and we consider the family of Cauchy problems 
 \begin{equation}
\label{e:nlcpr}
\left\{
\begin{array}{ll}
  \partial_t \rho_{\ee } + \partial_x \big[  \rho_{\ee } \ b( \rho_{\ee } \ast \eta_\ee) \big] = 0 \\
   \rho_{\ee } (0, x) = \bar  \rho(x). \\
\end{array}
\right.
\end{equation}
When $\ee \to 0^+$ (i.e. in the \emph{local} limit), the above Cauchy problem formally boils down to the conservation law
\begin{equation}
\label{e:clcpr}
\left\{
\begin{array}{ll}
  \partial_t \rho + \partial_x \big[ \rho \ b(\rho) \big] = 0  \\
  \rho (0, x) = \bar \rho(x). \\
\end{array}
\right.
\end{equation}
The by now classical theory by Kru{\v{z}}kov~\cite{Kruzkov} states that, if $\bar \rho \in L^\infty (\RR)$, the above problem has a unique \emph{entropy admissible} solution, i.e. 
loosely speaking a unique distributional solution that is consistent with the Second Principle of Thermodynamics.  
In~\cite{ACT}   P. Amorim, R. Colombo and A. Teixeira posed the following question. 
 \begin{question}
 \label{?} Can we rigorously justify the local limit? Namely, does the solution  $\rho_\ee$ of \eqref{e:nlcpr} converge to the entropy admissible solution $\rho$ of \eqref{e:clcpr} as $\ee \to 0^+$, in some suitable topology? 
\end{question}
In~\cite{ACT} the authors provide numerical evidence  supporting a positive answer to Question~\ref{?}. See also~\cite{ACG,BlandinGoatin,FKG,GoatinScialanga}. However, in~\cite{CCS} it is shown that the answer to Question~\ref{?} is, in general, negative.
More precisely, in~\cite{CCS} we exhibit some counterexamples ruling out convergence (see~\S\ref{s:ex} for an overview of these counterexamples). 

In~\cite{CCS} we also consider the ``viscous counterpart'' of Question~\ref{?}. More precisely, we fix a viscosity parameter $\nu>0$ and add a viscous second order term to the right hand side of both~\eqref{e:nlcpr} and~\eqref{e:clcpr}. We arrive at 
\begin{equation}
\label{e:vnlcpr}
\left\{
\begin{array}{ll}
   \partial_t \rho_{\ee \nu} + \partial_x \big[ \rho_{\ee \nu} b(\rho_{\ee \nu} \ast \eta_\ee) \big] = \nu \partial_{xx}^2 \rho_{\ee \nu} 
 \\
   \rho_{\ee \nu} (0, x) = \bar  \rho(x) \\
\end{array}
\right.
\end{equation} 
and  
\begin{equation}
\label{e:vclpr}
\left\{
\begin{array}{ll}
   \partial_t \rho_{\nu} + \partial_x \big[ \rho_{\nu} b(\rho_{ \nu}) \big] = \nu \partial_{xx}^2 \rho_{ \nu} 
 \\
   \rho_{\nu} (0, x) = \bar  \rho(x), \\
\end{array}
\right.
\end{equation} 
respectively. This yields the ``viscous counterpart'' of Question~\ref{?}, namely 
\begin{question}
 \label{?2} Fix $\nu>0$. Does the solution  $\rho_{\ee \nu}$ of \eqref{e:vnlcpr} converge to the solution $\rho_\nu$ of \eqref{e:vclpr}, when $\ee \to 0^+$? 
\end{question}
The answer to Question~\ref{?2} is largely positive. More precisely,~\cite[Theorem 1.1]{CCS} states in particular  that, for every $\nu>0$ and $T>0$, the family $\rho_{\ee \nu}$ converges to $\rho_\nu$ in the strong topology of $L^2 ([0, T] \times \RR)$\footnote{The precise results collected in~\cite[Theorem 1.1]{CCS} are actually stronger and in particular apply to the case of several space dimensions.}.  
To conclude the overview of the analytic results, we quote~\cite[Proposition~1.2]{CCS}, which establishes the ``nonlocal'' vanishing viscosity limit $\nu \to 0^+$ from~\eqref{e:vnlcpr} to~\eqref{e:nlcpr}, whereas  the ``local'' vanishing viscosity limit  from~\eqref{e:vclpr} to~\eqref{e:clcpr} is a classical result by Kru{\v{z}}kov~\cite{Kruzkov}. Summing up, we have the following convergence scheme: 
\begin{equation}
\label{e:disegno}
\minCDarrowwidth70pt
\begin{CD}
 \partial_t \rho_{\ee \nu} + \partial_x \big[ \rho_{\ee \nu} b(\rho_{\ee \nu} \ast \eta_\ee) \big] = \nu \partial_{xx}^2 \rho_{\ee \nu}     @> \ee \to 0^+  >  \text{\cite[Theorem 1.1]{CCS} }>  \partial_t \rho_\nu + \partial_x \big[ \rho_\nu b (\rho_\nu) \big] = \nu \partial_{xx}^2 \rho_\nu \\
@V \nu \to 0^+  V \text{\cite[Proposition~1.2]{CCS}} 
V        @V \nu \to 0^+ V \text{Kru{\v{z}}kov's Theorem} V\\
 \partial_t \rho_{\ee } + \partial_x \big[ \rho_{\ee } b(\rho_{\ee } \ast \eta_\ee) \big] = 0      @>  \ee \to 0^+ > \text{False in general}>   \partial_t \rho + \partial_x \big[ \rho b(\rho) \big] = 0.
\end{CD}
\end{equation}
\subsection{Numerical results} \label{ss:inr}
As pointed out before, the numerical evidence exhibited in~\cite{ACT} supports a positive answer to Question~\ref{?}, but this is contradicted by the analytic counterexamples in~\cite{CCS}. The present work aims at providing insights on the reason why the numerical evidence provides the wrong intuition. 

First, we point out that the numerical results in~\cite{ACT} are obtained by Lax-Friedrichs type schemes, which are known to have a very high numerical viscosity, see~\cite{Tadmor}. We refer to~\cite{LeVeque} for a more extended discussion, but, very loosely speaking, the numerical viscosity is a collection of finite difference terms that is the ``numerical counterpart'' of a viscous second order term like the one at the right hand side of the equations in~\eqref{e:vnlcpr} and~\eqref{e:vclpr}. In other words, the presence of the numerical viscosity implies that the model equation for the Lax-Friedrichs scheme applied to the conservation law at the first line of~\eqref{e:clcpr} is actually the equation at the first line of~\eqref{e:vclpr}, where the coefficient $\nu$ is of the same order as the space mesh.  Similarly, when the Lax-Friedrichs scheme is applied to the nonlocal conservation law at the first line of~\eqref{e:nlcpr} the model equation is actually the equation at the first line of~\eqref{e:vnlcpr}. 

We can now go back to the fact that the numerical evidence is not consistent with the analytic results: a possible explanation is the following. 
Because of the numerical viscosity, what the numerical tests in~\cite{ACT} are actually capturing is the convergence of $\rho_{\ee \nu}$ to $\rho_\nu$, which holds true by~\cite[Theorem 1.1]{CCS}. In other words: the numerical tests were designed to provide an answer to Question~\ref{?}, but as a matter of fact, owing to the numerical viscosity, they provide an answer to Question~\ref{?2}. Since the two questions have opposite answers, the numerical tests provide the wrong intuition concerning Question~\ref{?}. 

In the present paper we exhibit numerical experiments supporting the previous argument. In particular, we show that the numerical viscosity jeopardizes 
the reliability of standard numerical schemes for the study of the nonlocal-to-local limit from~\eqref{e:nlcpr} to~\eqref{e:clcpr}. In particular, in~\S\ref{ss:NumRes_Test1},~\S\ref{ss:NumRes_Test3},~\S\ref{ss:NumRes_Test2} we consider the 
counterexamples exhibited in~\cite{CCS} to show that the answer to Question~\ref{?} is \emph{negative} and we test them with the Lax-Friedrichs type scheme. The numerical results we obtain strongly suggest that the answer to Question~\ref{?} is \emph{positive} and hence provide the \emph{wrong} intuition. 

In this work we also further investigate the role of numerical viscosity  by comparing the Lax-Friedrichs type scheme with a Godunov type scheme. Lax-Friedrichs type schemes are known 
to have higher numerical viscosity then Godunov type schemes, see~\cite{Tadmor}. Consistently, we find that in several cases the numerical results obtained with the Godunov type scheme are better (i.e., more consistent with the analytic results) than those obtained with the Lax-Friedrichs type scheme, see again~\S\ref{ss:NumRes_Test1},~\S\ref{ss:NumRes_Test3},~\S\ref{ss:NumRes_Test2}. Finally, we provide further insights on the relation between numerical viscosity and nonlocal-to-local limit by varying the relation between the convolution parameter $\ee$ and the numerical viscosity.  Since the numerical viscosity depends monotonically on the space mesh, it suffices to vary the relation between the convolution parameter $\ee$ and the space mesh $h$. The numerical results obtained when $h$ is of the order of $\ee^2$ are better (i.e., more consistent with the analytic results) than those obtained 
when $h$ is of the order of $\ee$, see~\S\ref{ss:NumRes_Test1} and~\S\ref{ss:NumRes_Test2}. This again shows that the numerical viscosity compromises the reliability of the numerical schemes: indeed, when the numerical viscosity decays faster to $0$ the numerical results are more reliable.  In general, the best results are obtained with the Godunov type scheme when the space mesh $h$ is of the order of $\ee^2$, and this again confirms that the smaller the numerical viscosity, the more reliable the numerical results. 

The paper is organized as follows. In~\S\ref{s:Scheme} we discuss the numerical schemes used in the present work, i.e. the Lax-Friedrichs and the Godunov schemes. In~\S\ref{s:ex} we introduce the examples we will use in the numerical tests and we overview their main analytic properties.  In~\S\ref{s:benchmark} we validate our schemes by computing the numerical solutions in examples where the analytic solution is known, and by showing that the two are close. In~\S\ref{s:NumRes} we introduce our main numerical results concerning the nonlocal-to-local limit. In~\S\ref{s:end} we draw our conclusions and we outline some possible future work. 
To simplify the exposition, in the paper we always focus on the case where the conservation law at the first line of~\eqref{e:clcpr} is the scalar Burgers' equation
\begin{equation}
\label{e:burgers}
        \partial_t \rho + \partial_x ( \rho^2)  =0
\end{equation} 
and hence the nonlocal equation at the first line of~\eqref{e:nlcpr} is 
\begin{equation}
\label{e:nlburgers}
        \partial_t \rho_\ee + \partial_x \big( \rho_\ee (\rho_\ee \ast \eta_\ee ) \big)  =0. 
\end{equation}

\section{Two numerical schemes for the Burgers' equation}
\label{s:Scheme}
We now discuss two numerical schemes for both the local~\eqref{e:burgers} and nonlocal Burgers' equation~\eqref{e:nlburgers}. We refer to the book by LeVeque~\cite{LeVeque} for an extended discussion on numerical schemes  for conservation laws.

We discretize the $(t, x)$-plane by choosing the space mesh width $h$  and the time step $\Delta t$ and by introducing the 
mesh points $(t_n, x_j)$ given by $x_j = j h$, $j \in \mathbb Z$, and $t_n = n \Delta t$, $n=0, \dots, N$, $N= [T/ \Delta t] + 1$, where $T$ is the final time and $[ \cdot ]$ denotes the integer part. In the following we always consider a uniform mesh where $\Delta t / h= 1/6$, which is consistent with the CFL condition. 
For technical reasons we also define 
$$
    x_{j + 1/2} = x_j + h/2 = ( j + 1/2)h.   
$$ 
%Note that in our examples the modulus of the initial data will be always bounded by $1$ and hence, owing to the maximum principle,  
%the maximum propagation speed for the conservation law~\eqref{e:burgers} is bounded by $2$. This implies that the choice $\Delta t / h= 1/6$
%is consistent with the Courant-Friedrichs-Lewy (CFL) condition for the conservation law~\eqref{e:burgers}. The nonlocal equation~\eqref{e:nlburgers} does not have in general maximum principle and the a-priori bounds on the propagation speed are $\ee$-dependent and blow up as $\ee \to 0^+$. 
%This implies that the choice $\Delta t / h= 1/6$ might violate the CFL condition for~\eqref{e:nlburgers}. This issue might be the subject of further investigation. Since in the present work we want to compare our results with those in~\cite{ACT}, we use a similar uniform 
%mesh.
%{\color{magenta}Note by Laura: do you agree with the previous comment? Should we erase it?} {\color{cyan}I think we can keep this comment.}
The numerical schemes aim at defining a piecewise constant approximate solution $\rho_h$. As a matter of fact, in the following we will only define the discrete values $\rho^n_j$. The pointwise values of the 
approximate solutions are recovered by setting $\rho_h (t, x) : = \rho^n_j$ if $(t, x) \in ]t_{n}, t_{n+1} [ \times ] x_{j-1/2}, x_{j+ 1/2}[.$ We construct the approximate initial datum by setting 
$$
     \rho^0_j : =  \frac{1}{h} \int_{x_{j - 1/2}}^{x_{j + 1/2}} \bar \rho (x) dx.  
$$
Both the Lax-Friedrichs and the Godunov scheme are \emph{conservative} methods that can be written in the form 
\begin{equation}
\label{e:conservative}
        \rho^{n+1}_j = \rho_j^n- \frac{\Delta t}{h} \Big[ F^n_{j+1/2} - F^n_{j - 1/2}  \Big],
\end{equation}
where $F^n_{j+1/2}$ is the so-called numerical flux function. The two methods differ in the way one defines 
the value of $F^n_{j+1/2}$. We now separately describe them.  
\subsection{The Lax-Friedrichs method}
	\label{ss:Scheme_LF}
	The Lax-Friedrichs scheme was originally designed for the nonlinear conservation law
	%{\color{cyan}Note: in this sentence, we can understand that LF is designed especially for CL, but it can be used for multiple equations...}
	\begin{equation} 
	\label{e:cl}
	   \partial_t \rho + \partial_x f(\rho) =0  
	\end{equation}
	and it is 
	defined by plugging into~\eqref{e:conservative} the following numerical flux function:  
\begin{equation} \label{eq:LaxFriedrichs}
		F^n_{j+1/2} =  \frac{h}{2 \Delta t} (\rho^n_{j} - \rho^n_{j+1}) + \frac{1}{2} (f( \rho^n_{j}) + f(\rho^n_{j+1}) ).  
\end{equation}
In the case of the Burgers' equation~\eqref{e:burgers}, the above expression boils down to   
\begin{equation} \label{eq:LF_Dirac}
		F^n_{j+1/2} =  \frac{h}{2 \Delta t} (\rho^n_{j} - \rho^n_{j+1}) + \frac{1}{2} (( \rho^n_{j})^2 + (\rho^n_{j+1})^2 ) \implies 
		\rho_j^{n+1} \ = \ \dfrac{1}{2}(\rho_{j+1}^n + \rho_{j-1}^n) - \dfrac{\Delta t}{2h}\left[ (\rho_{j+1}^n)^2 - (\rho_{j-1}^n)^2 \right].
\end{equation}
Note that the Lax-Friedrichs scheme is first order accurate in both time and space.  Also, the numerical viscosity $\nu_{LF}$ satisfies $\nu_{LF}= h^2 / 2 \Delta t$, and owing to the CFL condition $\Delta t/h = 1/6$ we get $\nu_{LF}=3h$. In other words, the numerical viscosity  is of the same order as the same space mesh.

Lax-Friedrichs type schemes for nonlocal conservation laws were considered in various works, see for instance in~\cite{ACT,Sedimentation,BlandinGoatin}. In the case of the 
nonlocal Burgers' equation~\eqref{e:nlburgers}, the numerical flux function is defined by setting 
\begin{equation} \label{eq:nlLaxFriedrichs}
		F^n_{j+1/2} =  \frac{h}{2 \Delta t} (\rho^n_{j} - \rho^n_{j+1}) + \frac{1}{2} (\rho^n_{j} c^n_j + \rho^n_{j+1} c^n_{j+1} ),
\end{equation}
where $c_j^n $ is the approximate value of the convolution kernel and in the present work it is computed by the quadrature formula 
\begin{equation} \label{e:cienne}
    c_j^n = \sum_{k= - \ell}^{\ell-1} \gamma_k \rho^n_{j-k}, \quad 
    \text{where} \; 
    \gamma_k = \int_{kh}^{(k+1)h} \eta_\ee (y) dy \; \text{and} \; \ell = \left[ \frac{\ee}{h}\right] + 1. 
\end{equation}
We recall that the support of the convolution kernel $\eta_\ee$ is always contained in the interval $[- \ee, \ee]$. By plugging~\eqref{eq:nlLaxFriedrichs} 
into~\eqref{e:conservative} we arrive at    
\begin{equation} \label{eq:LF_Kernel}
		\rho_j^{n+1} \ = \ \dfrac{1}{2}(\rho_{j+1}^n + \rho_{j-1}^n) - \dfrac{\Delta t}{2h}\left( \rho_{j+1}^n c_{j+1}^n - \rho_{j-1}^n c_{j-1}^n \right). 
\end{equation}

\subsection{The Godunov method}
	\label{ss:Scheme_Godunov}
The basic idea underpinning the Godunov scheme is to solve Riemann problems on each cell of the computational mesh.  More precisely, the Godunov scheme for the nonlinear conservation law~\eqref{e:cl} is obtained by plugging the numerical flux function 
\begin{equation}
\label{e:ffGodunov}
       F^n_{j+1/2} : = 	f \big(\rho^\ast(\rho_j^n, \rho^n_{j+1}) \big)
\end{equation}
into~\eqref{e:conservative}. In the previous expression, $\rho^\ast (\rho_j^n, \rho^n_{j+1})$ is the value at the line $x=0$ of the entropy admissible solution of the Riemann problem between $\rho^n_j$ (on the left) and $\rho^n_{j+1}$ (on the right). Note that, owing to the Rankine-Hugoniot conditions, even if the solution of the Riemann problem has a discontinuity at $x=0$, the function $f(\rho)$ is continuous at $x=0$ and hence the value $f(\rho^\ast)$ is well-defined. As a matter of fact, if the flux function $f$ is convex we have the equality  
\begin{equation}
\label{e:rhoast}
       f \big(\rho^\ast(\rho_j^n, \rho^n_{j+1}) \big) =
      \left\{
    \begin{array}{ll}
    \min_{\rho \in [\rho_j^n, \rho_{j+1}^n]} f(\rho) & \rho_j^n \leq \rho^n_{j+1} \\
     \max_{\rho \in [\rho_{j+1}^n, \rho_{j}^n]} f(\rho)  & \rho_j^n \ge \rho^n_{j+1}, \\
    \end{array}
    \right.
\end{equation}
which in the case of the scalar Burgers' equation~\eqref{e:burgers} implies
\begin{equation}
\label{e:rhoastb}
   \big(\rho^\ast(\rho_j^n, \rho^n_{j+1}) \big)^2  =
      \left\{
    \begin{array}{ll}
    \min_{\rho \in [\rho_j^n, \rho_{j+1}^n]} \rho^2 & \rho_j^n \leq \rho^n_{j+1} \\
     \max_{\rho \in [\rho_{j+1}^n, \rho_{j}^n]} \rho^2  & \rho_j^n \ge \rho^n_{j+1}. \\
    \end{array}
    \right.
\end{equation}
Note furthermore that the Godunov scheme is known to have less numerical viscosity than the Lax-Friedrichs scheme, see~\cite{Tadmor}. 	

Godunov type schemes for nonlocal equations have been considered in~\cite{ChiarelloGoatin,FKG}. To define a Godunov scheme for the nonlocal Burgers' equation we first define the convolution term 
\begin{equation}
\label{e:V}
    V^n_{j + 1/2} = \sum_{k=-\ell}^{\ell-1} {\gamma}_k \,\rho_{j-k+1}^n, \quad 
    \text{with~}  \gamma_k \text{~as in~\eqref{e:cienne}}.
    %\text{where} \;  \tilde{\gamma}_k  \ = \ \int_{kh}^{(k+1)h}  \eta_\varepsilon(y) dy
    %\; \text{and} \; \ell = \left[ \frac{\ee}{h}\right] + 1. 
\end{equation}
By plugging the formula $f(\rho)=  V^n_{j + 1/2} \rho$ into~\eqref{e:rhoast} and recalling~\eqref{e:ffGodunov} we arrive at 
\begin{equation}
\label{e:nlGodunov}
        F^n_{j+1/2}  = 
%         \left\{
%    \begin{array}{ll}
%    \min_{\rho \in [\rho_j^n, \rho_{j+1}^n]} V^n_{j + 1/2} \rho & \rho_j^n \leq \rho^n_{j+1} \\
%     \max_{\rho \in [\rho_{j+1}^n, \rho_{j}^n]} V^n_{j + 1/2} \rho  & \rho_j^n \ge \rho^n_{j+1} \\
%    \end{array}
%    \right.
%    =  
  \left\{ \begin{array}{ll}
        V^n_{j + 1/2} \rho_j^n & V^n_{j + 1/2}  \ge 0 \\
    V^n_{j + 1/2} \rho_{j+1}^n & V^n_{j + 1/2}  < 0. \\
    \end{array}
    \right.
\end{equation}
By plugging the above numerical flux function into~\eqref{e:conservative} we obtain a Godunov type scheme for the nonlocal Burgers' equation~\eqref{e:nlburgers}. Note that our scheme is slightly different from the one in~\cite{ChiarelloGoatin,FKG} because in~\cite{ChiarelloGoatin,FKG} the authors focus on the case where $\bar \rho \ge 0$, which implies that $\rho_j^n \ge 0$ for every $n$ and $j$ and hence  that $  V^n_{j + 1/2} \ge 0$. This in turn implies that only the first case at the right hand side of~\eqref{e:nlGodunov} can occur. In the present work we consider cases where $\bar \rho$ attains negative values (see Example A in~\S\ref{sss:A}) and hence we use~\eqref{e:nlGodunov}. 
\section{Analytic results} \label{s:ex} In this paragraph we briefly discuss the main analytic properties of the examples we will use in our numerical tests. 
\subsection{Example A: odd initial datum, isotropic convolution kernels}
\label{sss:A}
Assume that 
\begin{equation}
\label{eq:init_Test1}
		\rho(0, x)= \bar{\rho}^{(A)}(x) : =  (x+2)\mathbbm{1}_{[-2,-1]}(x) + 
		\mathbbm{1}_{[-1,0]}(x) - \mathbbm{1}_{[0,1]}(x) + (x-2) \mathbbm{1}_{[1,2]}(x),
\end{equation}
where here and in the following $\mathbbm{1}_E$ denotes the characteristic function of the set $E$. 
The entropy admissible solution of the Cauchy problem for the (local) Burgers' equation~\eqref{e:burgers} is  
\begin{equation} \label{eq:sol_Test1}
	\rho^{(A)}(t,x) = \left\lbrace \begin{array}{clllll}
	\dfrac{x+2}{2t+1}, & x \in [-2,2t-1], & t \leq \dfrac12, & \text{~or~} x \in [-2,0], & t > \dfrac12, \\[0.8em]
									1, & x \in \, [2t-1,0], & t \leq \dfrac12,  \\[0.8em]
									-1, & x \in \, [0, 1-2t], & t \leq \dfrac12, \\[0.8em]
      \dfrac{x-2}{2t+1}, & x \in \, [1-2t,2], & t \leq \dfrac12, & \text{~or~} x \in \, [0,2], & t > \dfrac12\\
      [0.8em] 
      0 & \text{elsewhere.} 
								\end{array} \right.
\end{equation}
We now consider  the Cauchy problem obtained by coupling~\eqref{eq:init_Test1} with the nonlocal Burgers' equation~\eqref{e:nlburgers} and we term $\rho_\ee^{(A)}$ its solution. We assume furthermore that the convolution kernel $\eta$ is even, i.e. $\eta_\ee (x) = \eta_\ee (-x)$ for every $x$. The analysis  in~\cite[\S5.1]{CCS} states that, under these assumptions, the family  $\rho_\ee^{(A)}$ does not converge to the entropy admissible solution~\eqref{eq:sol_Test1} as $\ee \to 0^+$, not even weakly or up to subsequences.  We refer to~\cite{CCS} for the precise statements and the technical proof, but loosely speaking the very basic idea is the following. By using the fact that the initial datum $\bar \rho^{(A)}$ is odd and that the convolution kernel is even, one can show that the solution of the nonlocal equation is odd and this in turn implies, after some more work, that 
\begin{equation}\label{e:cons}
    \int_{-\infty}^0 \rho_\ee^{(A)}(t, y) dy = 
    \int_{-\infty}^0 \bar \rho^{(A)}( y) dy, \quad \text{for every $t>0, \ee>0$.}
\end{equation} 
On the other hand, the entropy admissible solution of the Burgers' equation satisfies 
\begin{equation}\label{e:ncons}
    \int_{-\infty}^0 \rho^{(A)}(t, y) dy < 
    \int_{-\infty}^0 \bar \rho^{(A)}( y) dy \quad \text{for every $t>0$}
\end{equation} 
and by comparing~\eqref{e:cons} and~\eqref{e:ncons} and performing some more work one eventually manages to rule out convergence.  

\subsection{Example B: positive initial datum, anisotropic convolution kernels} \label{exB}
If 
\begin{equation}
\label{eq:init_Test3}
		\rho(0, x)= \bar \rho^{(B)}(x) : =  \mathbbm{1}_{[-1,0]}(x),
\end{equation}
then the entropy admissible solution of the Cauchy problem for the (local) Burgers' equation~\eqref{e:burgers} 
is 
\begin{equation} \label{eq:sol_Test3}
	\rho^{(B)}(t,x) = \left\lbrace \begin{array}{cll}
									\dfrac{x+1}{2t}, & x \in [-1,2t-1], & t \leq 1 \\[0.8em]
									1, & x \in \, ]2t-1,t], & t \leq 1 \\[0.5em]
									\dfrac{x+1}{2t}, & x \in [-1,2\sqrt{t}-1], & t > 1\\[0.8em]
									0 & \text{elsewhere}.
								\end{array} \right.
\end{equation}
We term $\rho^{(B)}_\ee(t,x)$ the solution of the Cauchy problem obtained by coupling~\eqref{e:nlburgers} with~\eqref{eq:init_Test3}.  Assume that the convolution kernels $\eta_\ee$ are \emph{anisotropic}, more precisely they are supported on the negative real line, i.e. 
\begin{equation} \label{e:support} 
     \eta_\ee(x) =0 \quad \text{for every $x>0$}.
\end{equation}     
In this case the analysis  in~\cite[\S5.2]{CCS} states that the family  $\rho_\ee^{(B)}$ does not converge to the entropy 
 admissible solution~\eqref{eq:sol_Test2} as $\ee \to 0^+$, not even weakly or up to subsequences. The 
basic reason why $\rho_\ee^{(B)}$ does not converge to $\rho^{(B)}$ is because one can show that 
\begin{equation} \label{e:rhoci}
\rho_\ee^{(B)}(t, x)=0 \quad \text{for every $x>0$ and $t>0$},
\end{equation}
see~\cite[Lemma 5.3]{CCS}. Since $\rho^{(B)}$ does not share this property, then with some more work one manages to rule out convergence.

\subsection{Example C: positive initial datum, isotropic convolution kernels} \label{exC}
If 
\begin{equation}
\label{eq:init_Test2}
		\rho(0, x)= \bar \rho^{(C)}(x) : =  \mathbbm{1}_{[-1,1]}(x),
\end{equation}
then the entropy admissible solution of the Cauchy problem for the (local) Burgers' equation~\eqref{e:burgers} 
is 
\begin{equation} 
\label{eq:sol_Test2}
	\rho^{(C)}(t,x) = \left\lbrace \begin{array}{cll}
									\dfrac{x+1}{2t}, & x \in [-1,2t-1], & t \leq 2, \\[0.8em]
									1, & x \in \, [2t-1,t+1], & t \leq 2, \\[0.8em]
									\dfrac{x+1}{2t}, & x \in [-1,2\sqrt{2t}-1], & t > 2 \\[0.8em]
									0 & \text{elsewhere.} 
								\end{array} \right.
\end{equation}
As before, we term $\rho^{(C)}_\ee(t,x) $ the solution of the Cauchy problem obtained by coupling~\eqref{e:nlburgers} with~\eqref{eq:init_Test2}.
 Assume that the convolution kernels are  even functions, i.e. $\eta_\ee (x) = \eta_\ee (-x)$, for every $x \in \mathbb{R}$. In this case, for every $p>1$ the analysis in~\cite[\S5.2]{CCS} 
states that, as $\ee \to 0^+$, $\rho^{(C)}_\ee$ does not converge to $\rho^{(C)}$ strongly in $L^{p}$, not even up to subsequences.  
Loosely speaking, this is due to the fact that we can single out an entropy that is conserved by $\rho^{(C)}_\ee$ and is dissipated by $\rho^{(C)}$.

\subsection{Example D: explicit solution of the nonlocal equations} \label{sss:D}
If 
\begin{equation}
\label{eq:init_TestD}
		\rho(0, x)= \bar \rho^{(D)}(x) : =  \mathbbm{1}_{]-\infty,0]}(x),
\end{equation}
then the entropy admissible solution of the Cauchy problem for the (local) Burgers' equation~\eqref{e:burgers} 
is the shock 
\begin{equation} \label{eq:sol_TestD}
	\rho^{(D)}(t,x) = 
	\left\lbrace \begin{array}{ll} 1 & x  \in ] - \infty, t] \\
	     0 & x \in [t, + \infty[ . \\
								\end{array} \right.
	\end{equation}
Also, consider the nonlocal Burgers' equation~\eqref{e:nlburgers} and assume that the convolution kernel is supported on the positive real axis, i.e. $\eta_\ee (x)=0$ for every $x<0$.  In this case one can  show that, for every $\ee>0$, the solution of the  Cauchy problem obtained by coupling~\eqref{e:nlburgers} with~\eqref{eq:init_TestD} is exactly the same shock as in~\eqref{eq:sol_TestD}, i.e. $\rho^{(D)}_\ee \equiv \rho^{(D)}$. 

%We remark in passing that by introducing a suitable change of variables one can show that the case where $\eta_\ee (x) =0$ for every $x<0$ corresponds, in traffic models, to the case where drivers decide their speed based only on the downstream traffic density. Very loosely speaking, the reason why we have to introduce a change of variables is that the flux function of the (local) Burgers' equation is convex, whereas in (local) traffic models the flux function is concave. This case is very relevant in view of applications, has been extensively studied in several works (see e.g.~\cite{BlandinGoatin,ChiarelloGoatin,FKG,GoatinScialanga}) and has better analytic properties than the general case, in particular solutions of the nonlocal equations satisfy a maximum principle. To the best of our knowledge, the problem of the local limit $\ee \to 0^+$ is in this case open. Note, however, that a recent counterexample rules out the most ``natural'' convergence proof, see~\cite{CCS2}. 

\subsection{Example E: isotropic convolution kernels, regular limit solution} \label{sss:E} 
Assume that 
\begin{equation}
\label{eq:init_TestE}
		\rho(0, x)= \bar \rho^{(E)}(x) : = \frac{1}{4}\left(1+\sin\left(\frac{\pi x}{2}+\frac{\pi}{2}\right)\right)\mathbbm{1}_{[-2,0]}(x) + \frac12\mathbbm{1}_{[0,\infty[}(x).
\end{equation}
Since the initial datum is regular and monotone nondecreasing, classical results on scalar conservation laws rule out shock formation and imply that the solution of the Cauchy problem for the (local) Burgers' equation~\eqref{e:burgers} is regular.  Consider the nonlocal Burgers' equation~\eqref{e:nlburgers} and assume that the convolution kernels are even, i.e. that $\eta_\ee(x) = \eta_\ee (-x)$ for every $x \in \mathbb R$ and $\ee>0$. Owing to a convergence result by Zumbrun~\cite[Proposition 4.1]{Zumbrun}, in this case we expect that the solutions of the nonlocal equation uniformly converge to the solution of the (local) Burgers' equation.

\section{Benchmark numerical tests} \label{s:benchmark}
In this paragraph we discuss some benchmark tests we use to validate our numerical schemes. 
\subsection{Test 1: convergence of the numerical schemes for the local equation}
	\label{ss:Scheme_convergence}
\begin{figure}[!h]
	\centering \vspace*{-1.em}
	\includegraphics[width=\textwidth]{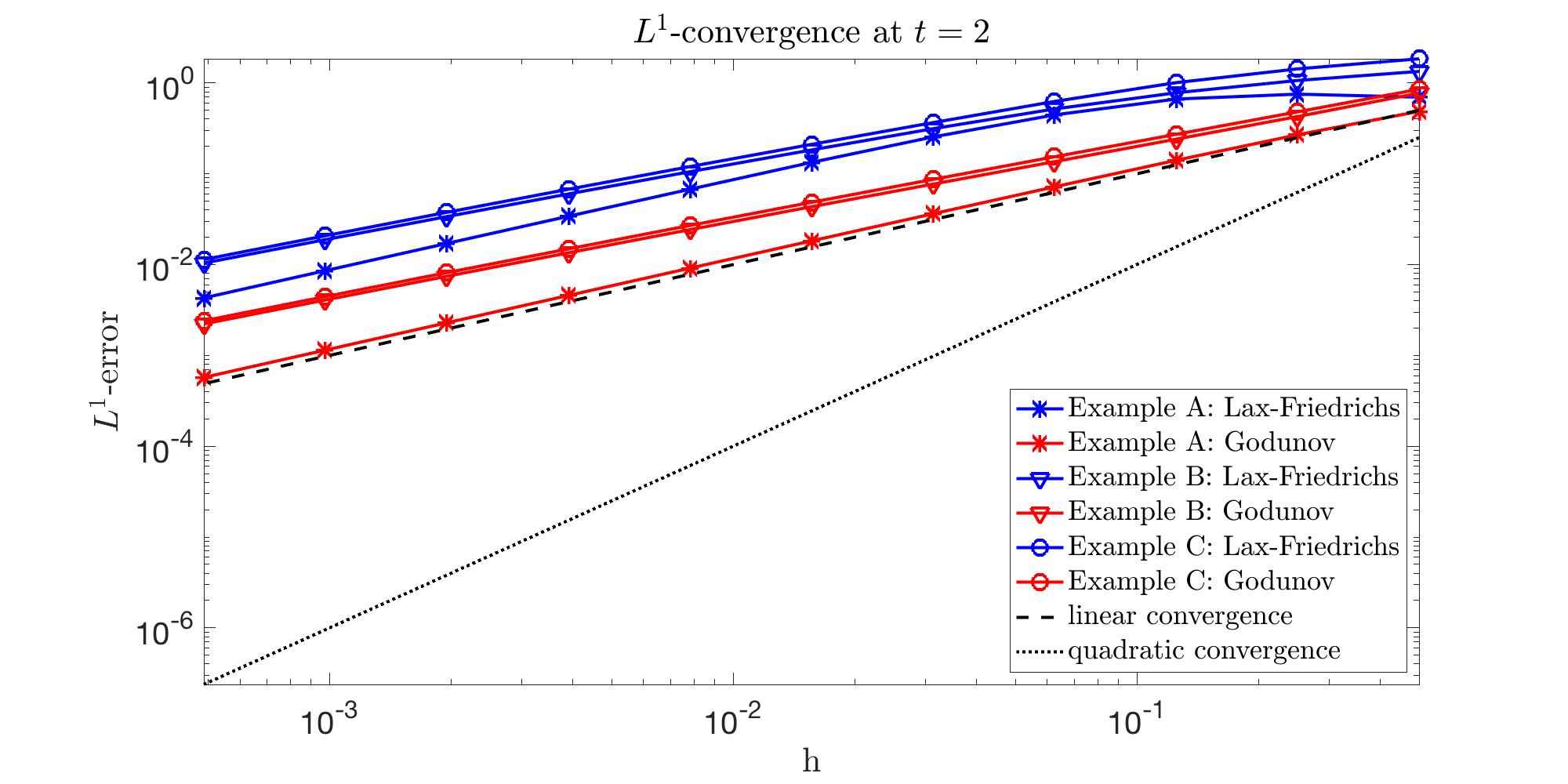} \vspace*{-1.em}
	\caption{\footnotesize Test 1, $L^1$-Convergence of the numerical schemes, Lax-Friedrichs and Godunov, with respect to the mesh size $h$, for the solution of the (local) Burgers' equation at $t = 2$ \label{fig:Scheme_convergence}}
\end{figure}

In Test 1, we focus on the (local) Burgers' equation~\eqref{e:burgers} and on same initial data $\bar \rho^{(A)}$, $\bar \rho^{(B)}$ and $\bar \rho^{(C)}$ as in~\eqref{eq:init_Test1},~\eqref{eq:init_Test3} and~\eqref{eq:init_Test2}, respectively.  We then compare the numerical solution given by the Lax-Friedrichs and Godunov schemes with the exact analytic entropy admissible solution, which in these cases can be explicitly computed, see~\eqref{eq:sol_Test1},~\eqref{eq:sol_Test3} and~\eqref{eq:sol_Test2}.  More precisely, we evaluate the $L^1$-norm in space of the difference between the exact and the numerical solution  
at fixed time $t=2$ and for different values of the space mesh. The results are displayed in 
 Figure~\ref{fig:Scheme_convergence} and show  the expected first order convergence for both schemes. 
	\subsection{Test 2 (Example D): validation of numerical schemes for the local and the nonlocal Burgers' equation}
	\label{ss:Other_Test6}
	Test 2 is designed to validate the numerical schemes for the nonlocal equation. 
	We take the same initial datum $\bar \rho^{(D)}$ as in~\eqref{eq:init_TestD} and the convolution kernel 
	\begin{equation} \label{eq:init_Test6}
	\begin{array}{rcl}
%		f(t,x,\rho) = \rho, & \qquad & v(r) = r,  \\[0.8em]
		\eta_\ee(x) : = \alpha_\ee (|x-\ee||x|)^{5/2}\mathbbm{1}_{[0,\ee]}(x), % , & & \rho^{(4)}_0(x) = \chi_{[-\infty,0]}(x). 
	\end{array}
\end{equation}
where (here and in the following) the constant $\alpha_\ee>0$ is chosen in such a way that $\eta_\ee$ has unit integral. The exact value of $\alpha_\ee$ can vary from occurrence to occurrence.  As pointed out in~\S\ref{sss:D}, in this case the solution of the Cauchy problem for the nonlocal Burgers' equation~\eqref{e:nlburgers} is explicit and it is given by~\eqref{eq:sol_TestD}. 	We can then validate the schemes for the nonlocal equation by computing the $L^1$ norm in space of the difference between the numerical solution (obtained with the Lax-Friedrichs  and the Godunov type schemes) and the exact analytic solution. 
\begin{figure}[!h]
	\centering
	\includegraphics[width=\textwidth]{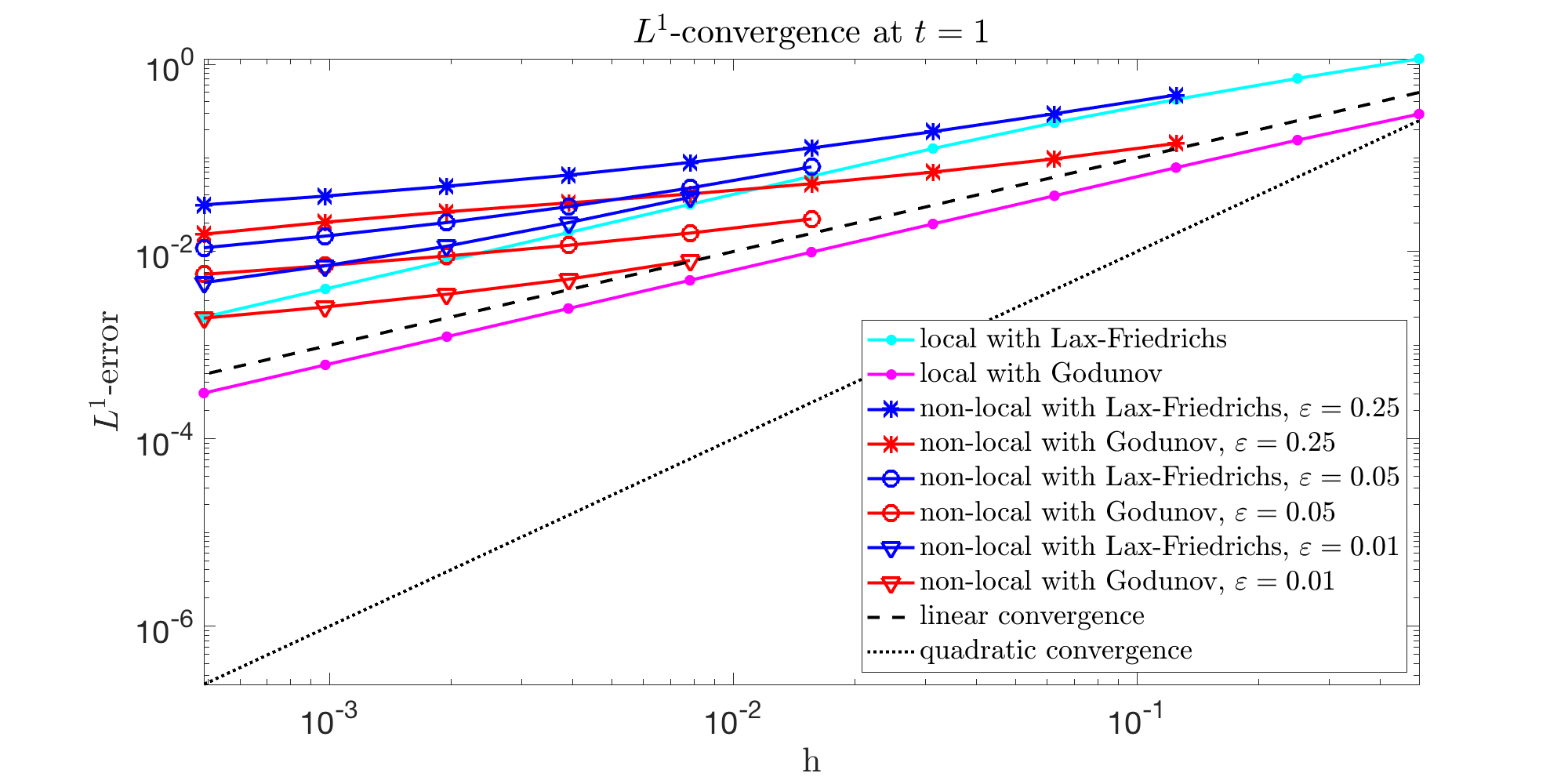}. 
	\caption{\footnotesize Test 2, $L^1$-Convergence of the numerical schemes, Lax-Friedrichs and Godunov, with respect to the mesh size $h$, for the (local) Burgers' equation and for the nonlocal equation at $\ee = 0.25, 0.05, 0.01$ and $t = 1$ \label{fig:Test6_conv-h}}
\end{figure}
The results are displayed in Figure~\ref{fig:Test6_conv-h}. We evaluate the $L^1$ norm at time $t=1$ for different values of the nonlocal parameter $\ee$, and we evaluate the convergence with respect to mesh size $h$. The results show convergence of the numerical solution to the analytic solution. 

\section{Numerical tests on the nonlocal-to-local limit}
\label{s:NumRes}
	\subsection{Test 3 (Example A): odd initial datum, isotropic convolution kernels}
	\label{ss:NumRes_Test1}
	In Test 3 we take the same initial datum $\bar \rho^{(A)}$ as in~\eqref{eq:init_Test1} and the convolution kernels   
\begin{equation}\label{e:even}
	 \eta_\ee(x) = \alpha_\ee (|x-\ee||x+\ee|)^{5/2}\mathbbm{1}_{[-\ee,\ee]}(x).
\end{equation}
As pointed out in~\S\ref{sss:A}, in this case the analysis in~\cite{CCS} implies that the solutions of the nonlocal equation~\eqref{e:nlburgers} \emph{do not} converge to the entropy admissible solution of the Burgers' equation, which is given by~\eqref{eq:sol_Test1}.  

\begin{figure}[!h]
	\centering
	\includegraphics[width=0.48\textwidth]{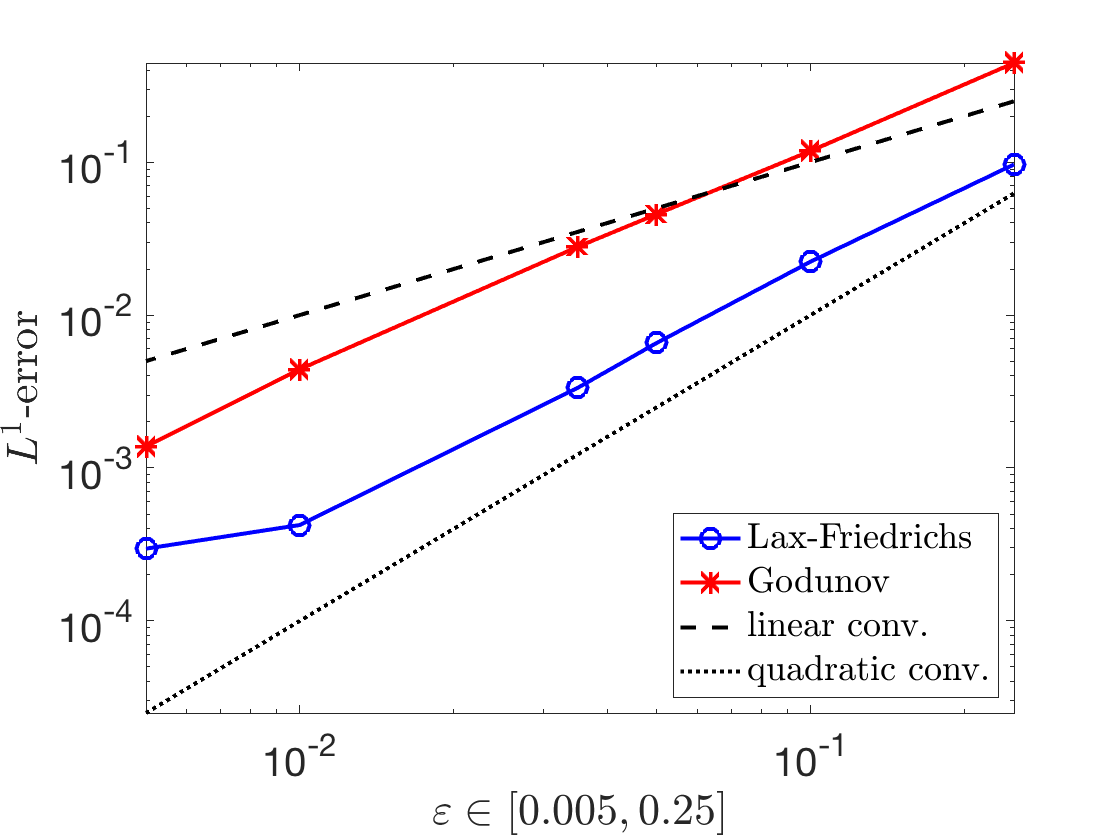} \hfill 
	\includegraphics[width=0.48\textwidth]{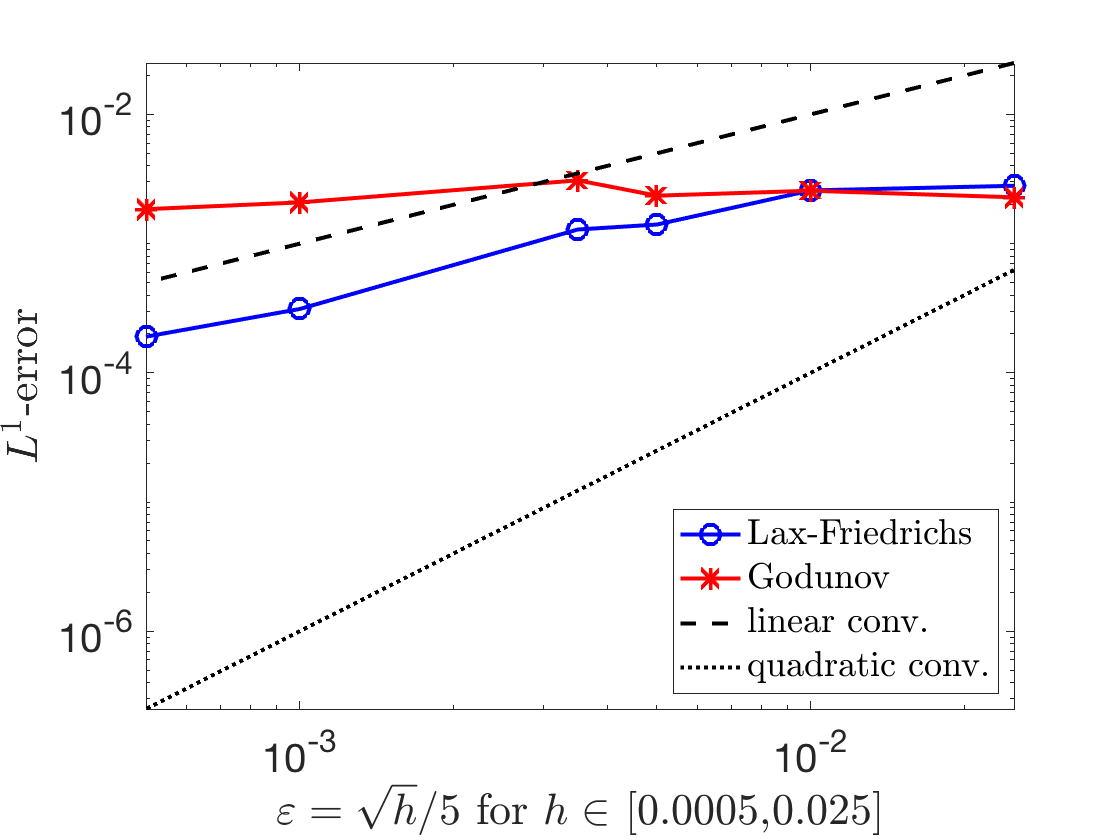}\\
	\hspace*{0.2\textwidth} (a) \hfill (b) \hspace*{0.2\textwidth} \textcolor{white}{.} 
	\caption{\footnotesize Test 3 (Example A), $L^1$-error at $t=2$, for different values of $\ee$, comparing the nonlocal solution to the local solution for both numerical schemes: (a) fixed viscosity $h=0.001$, (b) varying viscosity $h=25 \ee^2$. \label{fig:Test1_conv-a} }
\end{figure}
\begin{figure}[!h]
	\centering 
	\footnotesize{Analytic solution of the (local) Burgers' equation} \\
	\includegraphics[width=\textwidth]{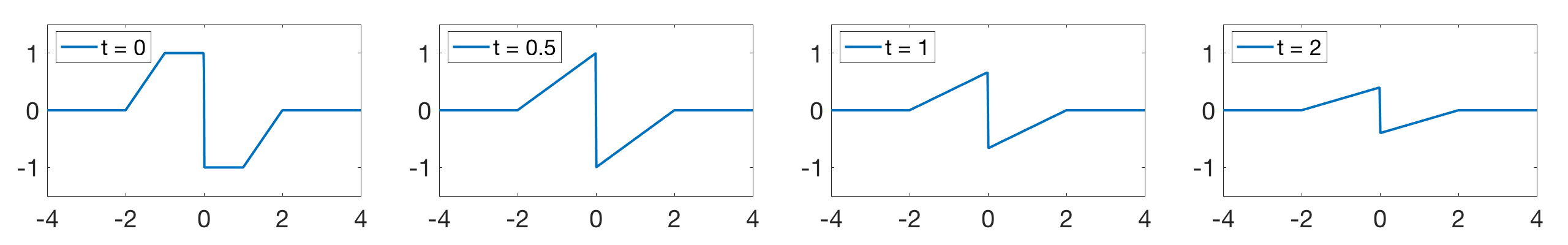} \\[1.em]
	\footnotesize{Numerical solution of the (local) Burgers' equation computed with the Lax-Friedrichs method} \\
	\includegraphics[width=\textwidth]{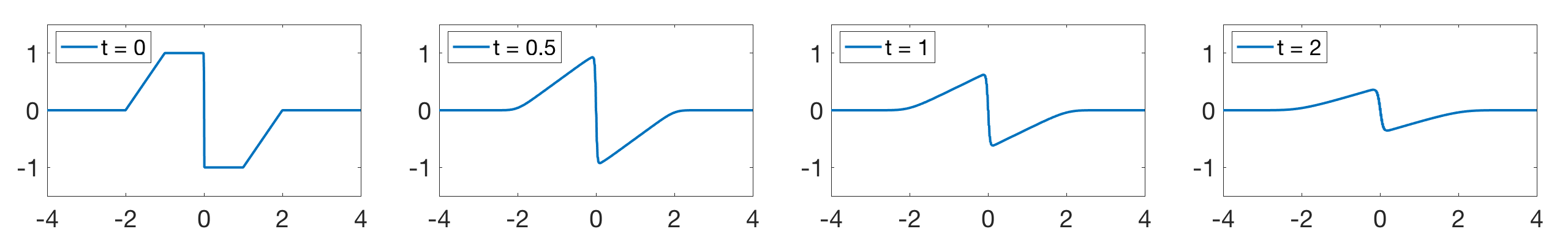} \\[1.em]
	\footnotesize{Numerical solution of the (local) Burgers' equation computed with the Godunov method} \\
	\includegraphics[width=\textwidth]{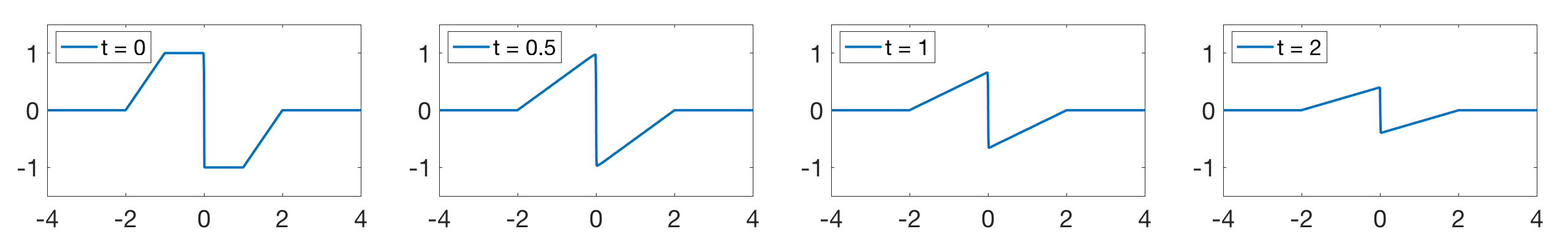} \\[1.em]
	\footnotesize{Numerical solution of the nonlocal Burgers' equation computed with the Lax-Friedrichs type method} \\
	\includegraphics[width=\textwidth]{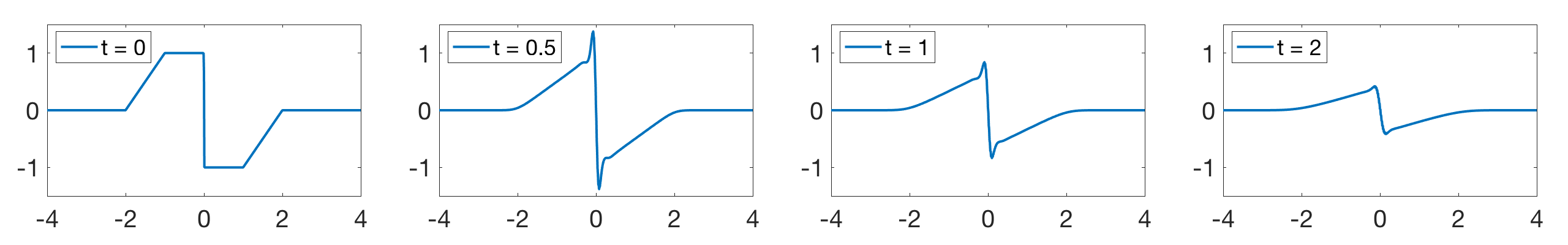} \\[1.em]
	\footnotesize{Numerical solution of the nonlocal Burgers' equation computed with the Godunov type method} \\
	\includegraphics[width=\textwidth]{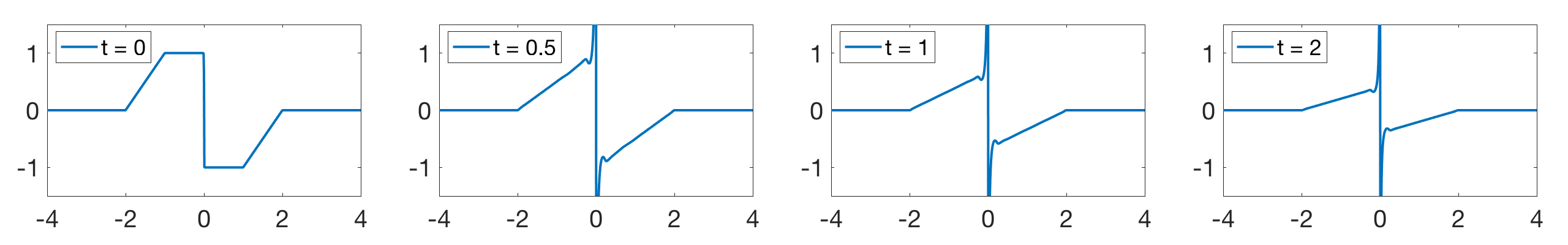}
	\caption{\footnotesize Test 3 (Example A), snapshots of solution of Burgers' equation with initial condition~\eqref{eq:init_Test1} and isotropic convolution kernel, when $\ee = 0.25$, $h = 0.01$. \label{fig:Test1}}
\end{figure}
In Test 3 we compute the numerical solution of the nonlocal equation by using the Lax-Friedrichs and the Godunov type schemes. Several snapshots of the solution are displayed in Figure~\ref{fig:Test1}.  
Also, we compare the numerical solution of the nonlocal equation with the analytic solution of the (local) Burgers' equation. More precisely, we evaluate the $L^1$ norm of the difference at time $t=2$ and at different values of the convolution parameter $\ee$.  We show the corresponding  results in Figure~\ref{fig:Test1_conv-a}. 

Here are the main remarks concerning the numerical results for Test 3.
\begin{itemize}
\item[i)] Figure~\ref{fig:Test1_conv-a}, part (a) shows the numerical results obtained by keeping the space mesh fixed and varying the 
convolution parameter $\ee$. The numerical results strongly suggest that the $L^1$ norm of the difference converges to $0$ when $\ee \to 0^+$. 
As pointed out in~\S\ref{sss:A}, in this case we can analytically rule out convergence and hence the numerical evidence provides the wrong intuition. Owing to the discussion in~\S\ref{ss:inr}, this is most likely    due to the presence of the numerical viscosity. 
%\item[ii)] Note furthermore that in Figure~\ref{fig:Test1_conv-a}, part (a), the convergence rate shown by the Godunov-type scheme is definitely slower than the one suggested by the Lax-Friedrichs type scheme, and this is consistent with the fact that Godunov schemes have lower viscosity than Lax-Friedrichs schemes, see~\cite{Tadmor}. 
\item[ii)]In Figure~\ref{fig:Test1_conv-a}, part (b), we display the results obtained by simultaneously varying the space mesh $h$ and the convolution parameter $\ee$. More precisely, we choose $h$ of the order of $\varepsilon^2$: in this way, the space mesh goes to $0$ much faster than the convolution parameter. The Lax-Friedrichs type scheme still suggests convergence, albeit at a slower rate than in Figure~\ref{fig:Test1_conv-a}, part (a). 
Conversely, the Godunov type scheme does not suggest convergence since the $L^1$ error is basically constant.  
These results are consistent with the fact that the numerical viscosity depends on the space mesh and it is higher in the Lax-Friedrichs scheme. Since the space mesh $h$ goes to $0$ very fast, the numerical schemes have less numerical viscosity and hence they provide a better intuition of the nonlocal-to-local limit.  Also, the Godunov scheme has lower numerical viscosity and hence it is more reliable than the Lax-Friedrichs scheme in investigating the nonlocal-to-local limit.
%\item[iii)]To conclude we comment on  the snapshots displayed in Figure~\ref{fig:Test1}. First, since the initial datum is odd, i.e. $\rho^{(A)}(x) = - \rho^{(A)}(-x)$, then the solutions of the nonlocal equation (viscous and inviscid) and of the local Burgers' equation (viscous and inviscid) are all odd, and this is consistent with the results in Figure~\ref{fig:Test1}. Also, the numerical solution of the nonlocal equation obtained with the Godunov type scheme has two high peaks at the left and at the right of $x=0$. This is consistent with the fact that by relying on analytic considerations  we expect high concentration of mass near the origin.
\end{itemize}
Wrapping up, Test 3 shows that the numerical viscosity jeopardizes the reliability of the numerical investigation of the nonlocal-to-local limit. To obtain reliable numerical results we had to introduce very low numerical viscosity. 
%%%%%%%%%%%%%%%%%%%
\subsection{Test 4 (Example B): positive initial datum, anisotropic convolution kernels}
	\label{ss:NumRes_Test3}	
	\begin{figure}[!h]
	\centering
	\includegraphics[width=\textwidth]{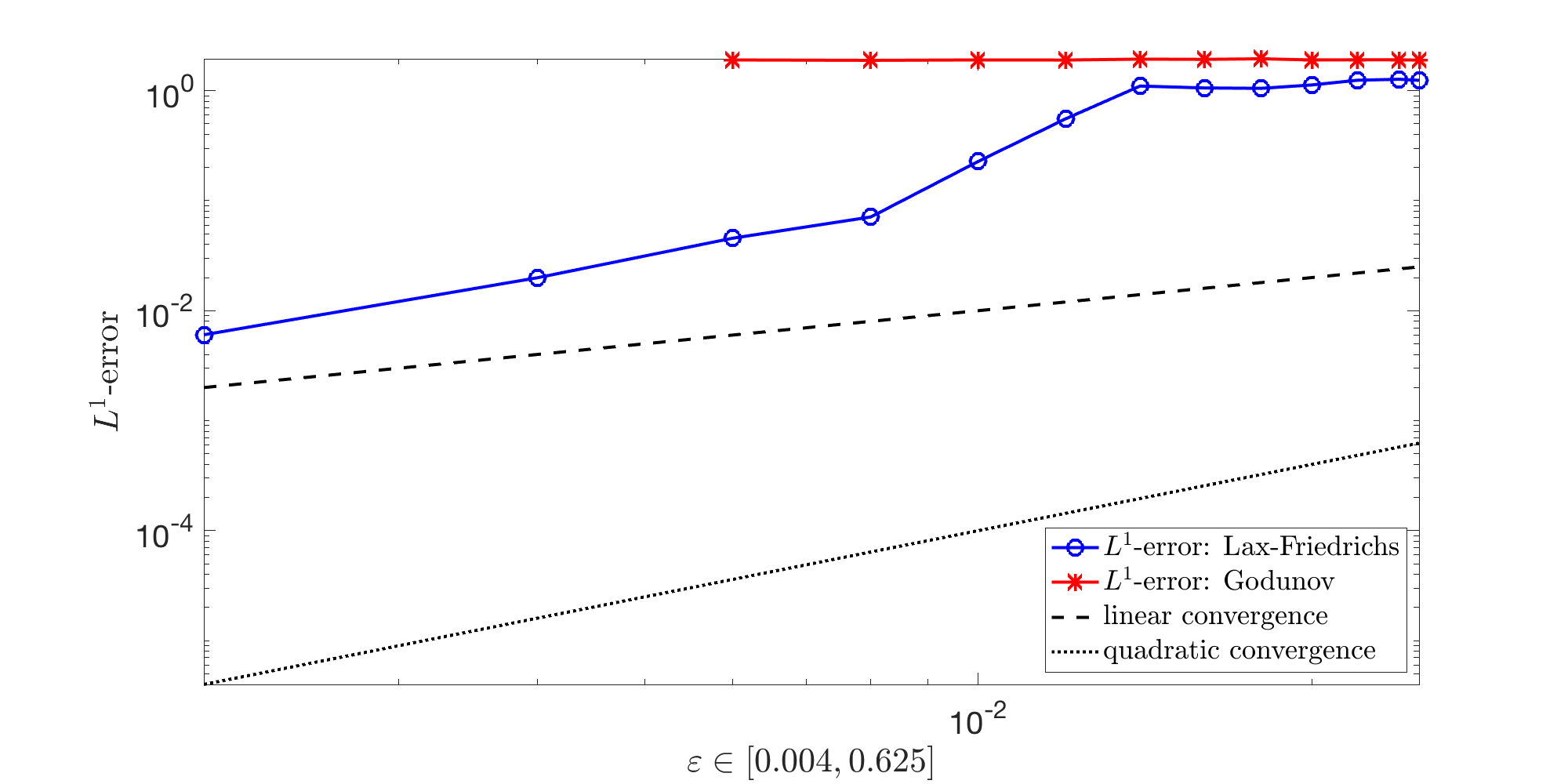}
	\caption{\footnotesize Test 4 (Example B),  $L^1$-error at $t=2$, for different values of $\ee$, comparing the solutions of the nonlocal equations with the entropy solution of the local equation for both numerical schemes and for varying viscosity $h$ such that $\ee = 1000 h^2$. The $L^1$ error of the Godunov scheme is much higher for small values of $\ee$.  \label{fig:Test3_conv-a-h} }
\end{figure}
	\begin{figure}[!h]
	\centering 
	\footnotesize{Analytic solution} \\
	\includegraphics[width=\textwidth]{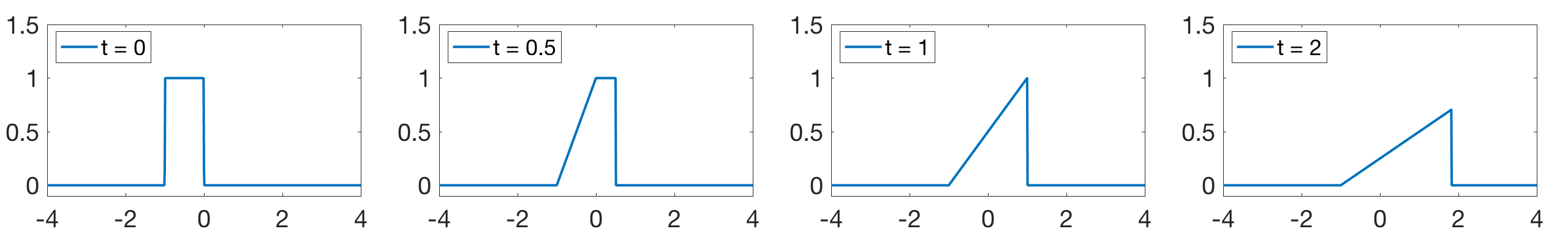} \\[1.em]
	\footnotesize{Numerical solution of the (local) Burgers' equation with Lax-Friedrichs method} \\
	\includegraphics[width=\textwidth]{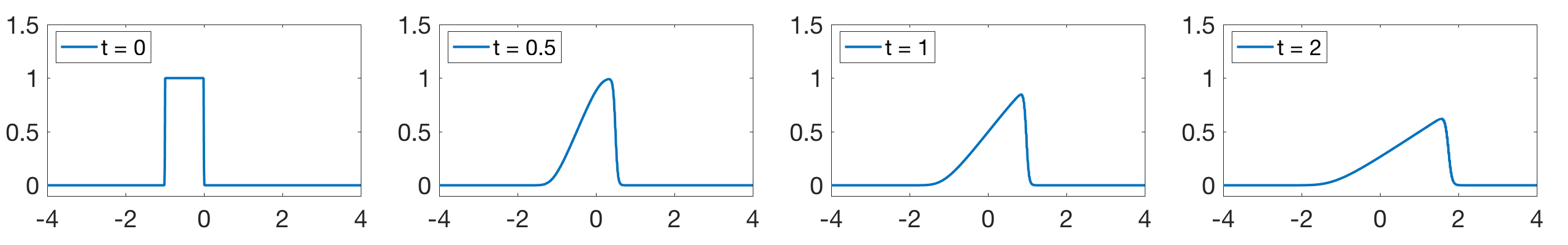} \\[1.em]
	\footnotesize{Numerical solution of the (local) Burgers' equation with Godunov method} \\
	\includegraphics[width=\textwidth]{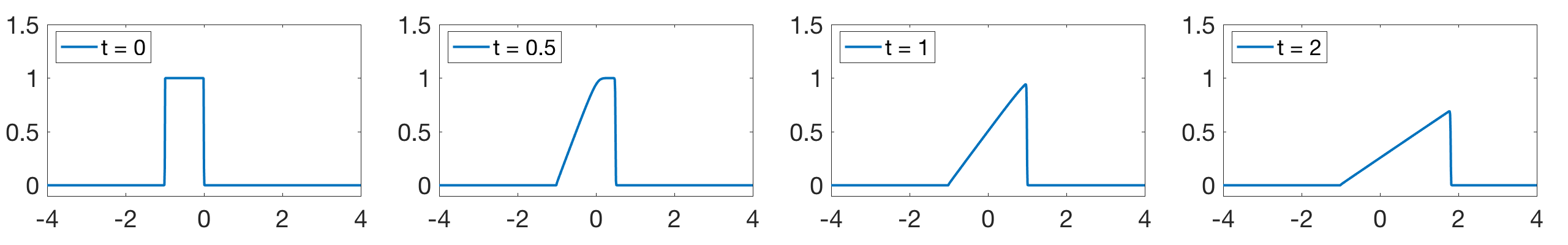} \\[1.em]
	\footnotesize{Numerical solution of the nonlocal Burgers' equation with Lax-Friedrichs method} \\
	\includegraphics[width=\textwidth]{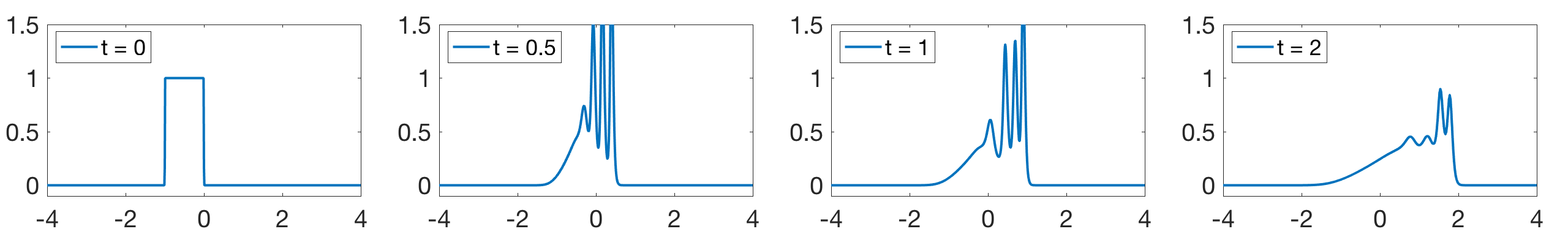} \\[1.em]
	\footnotesize{Numerical solution of the nonlocal Burgers' equation with Godunov method} \\
	\includegraphics[width=\textwidth]{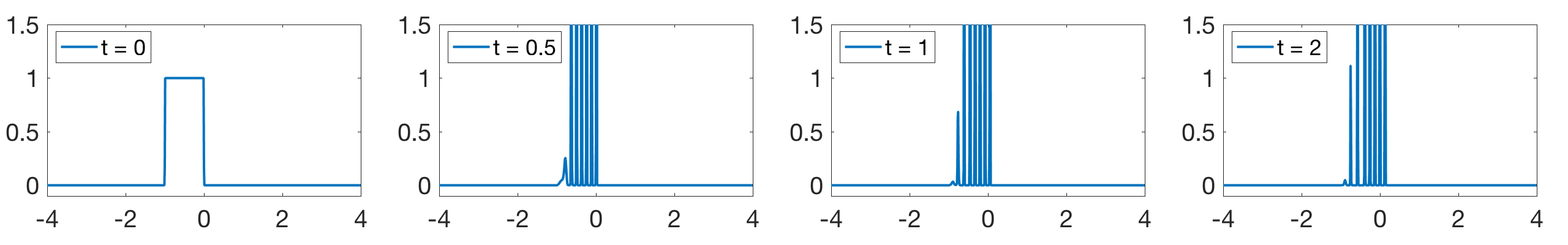} \\[1.em]
	\caption{\footnotesize Test 4 (Example B), snapshots of the solution of Burgers' equation with initial condition and convolution kernel~\eqref{eq:init_Test3}, when $\ee = 0.1$, $h = 0.01$. \label{fig:Test3}}
\end{figure}
In Test 4 we take the same initial datum $\bar \rho^{(B)}$ as in~\eqref{eq:init_Test3} and the convolution kernels
\begin{equation} \label{e:anisotropic} 	
\eta_\ee(x) = \alpha_\ee (|x||x+\ee|)^{5/2}\mathbbm{1}_{[-\ee,0]}(x). 
\end{equation}
Note that these convolution kernels satisfy~\eqref{e:support} and hence the discussion in~\S\ref{exB} applies. In particular, the analytic solutions $\rho^{(B)}_\ee$ of the nonlocal equations~\eqref{e:nlburgers} are all supported on the negative axis, i.e. satisfy~\eqref{e:rhoci}, and \emph{do not} converge to the solution of the (local) Burgers' equation. 

In Test 4 we compute the numerical solution of the nonlocal equations~\eqref{e:nlburgers} by using the  Lax-Friedrichs type method and the Godunov type method and we display the corresponding results in~Figure~\ref{fig:Test3_conv-a-h} and Figure~\ref{fig:Test3}. More precisely, Figure~\ref{fig:Test3_conv-a-h} displays the behavior of the $L^1$ norm of the difference between the (numerical) solutions of the nonlocal equation and the exact entropy admissible solution of the (local) Burgers' equation  at time $t=2$. Recall that the exact solution is given by~\eqref{eq:sol_Test3}. Figure~\ref{fig:Test3} shows the snapshots of the solution at time $t=2$. 
We now comment on~Figure~\ref{fig:Test3_conv-a-h} and Figure~\ref{fig:Test3}.
\begin{itemize}
\item[i)] The results in Figure~\ref{fig:Test3_conv-a-h} obtained with the  Lax-Friedrichs type method  suggest that the solutions of the nonlocal equation converge to  the entropy admissible solution of the Burgers' equation. This contradicts the analytic results discussed in~\S\ref{exB}. On the other hand, the numerical results obtained with the Godunov type scheme do not suggest convergence and hence 
are consistent with the analytic results in~\S\ref{exB}.  This is most likely  due to the fact that the Lax-Friedrichs type scheme has higher numerical viscosity and hence it is not reliable to test the nonlocal-to-local limit. The Godunov type scheme has less numerical viscosity and is therefore more reliable. 
\item[ii)] The snapshots of the solution obtained with the Godunov type and the Lax-Friedrichs type schemes confirm that the Godunov type scheme is more reliable. Indeed, the exact solution of the nonlocal Burgers' equation is supported on the negative axis, i.e. satisfies~\eqref{e:rhoci}.  Remarkably, this important analytic property is satisfied by the numerical solution obtained by the Godunov type method, but it is \emph{not} satisfied by the solutions obtained by the  Lax-Friedrichs type method, see Figure~\ref{fig:Test3}.  
\end{itemize}
The take-home message from Test 4 is the following: the Godunov type scheme is more reliable than the Lax-Friedrichs type scheme for the numerical investigation of the nonlocal-to-local limit. This is due to the fact that the Godunov type scheme is less affected by numerical viscosity, see~\cite{Tadmor}. 
\subsection{Test 5: smoother positive density and anisotropic convolution kernels}
\label{ss:Other_Test4}
\begin{figure}[!h]
	\centering 
	\footnotesize{Numerical solution of the (local) Burgers' equation computed with Lax-Friedrichs method} \\
	\includegraphics[width=\textwidth]{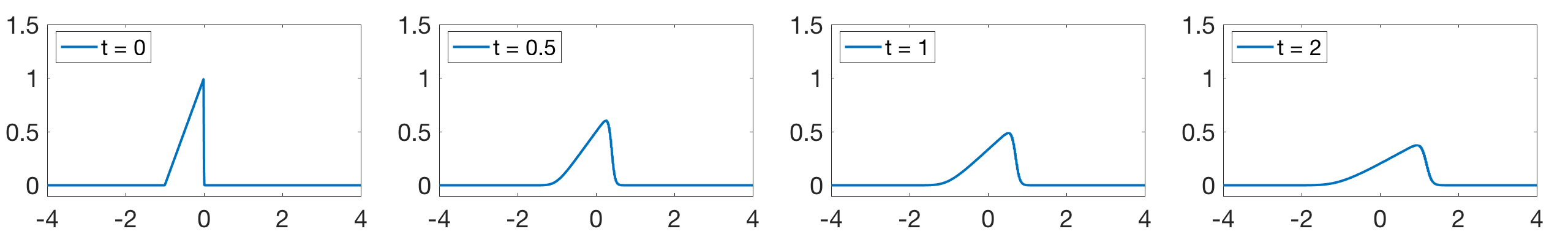} \\[1.em]
	\footnotesize{Numerical solution of the (local) Burgers' equation computed with Godunov method} \\
	\includegraphics[width=\textwidth]{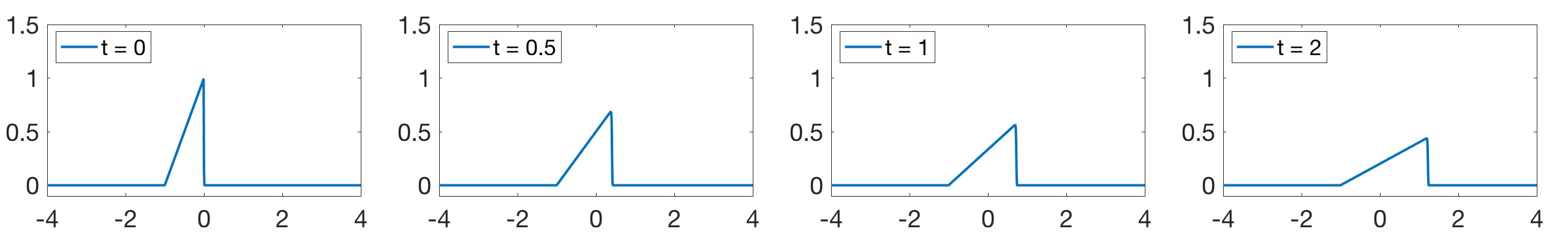} \\[1.em]
	\footnotesize{Numerical solution of the nonlocal Burgers' equation computed with Lax-Friedrichs method} \\
	\includegraphics[width=\textwidth]{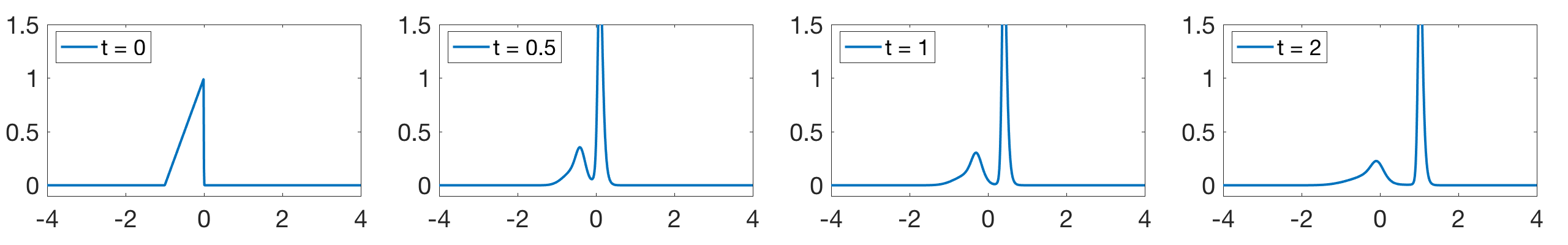} \\[1.em]
	\footnotesize{Numerical solution of the nonlocal Burgers' equation computed with Godunov method} \\
	\includegraphics[width=\textwidth]{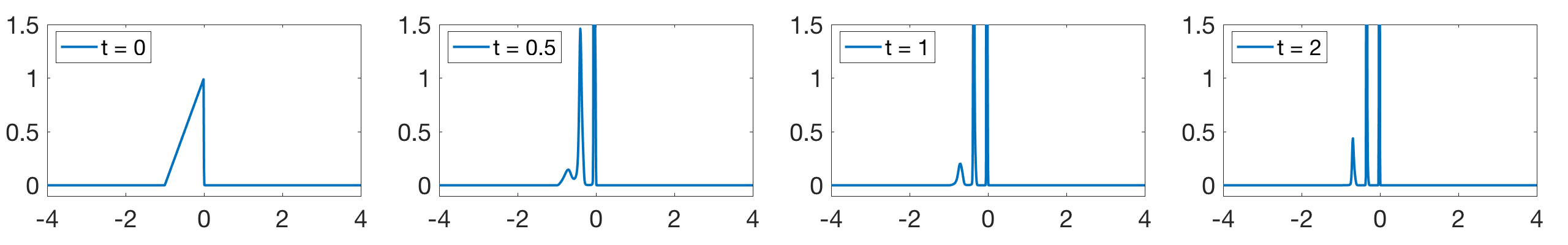} \\[1.em]
	\caption{\footnotesize Test 5, snapshots of solution of Burgers' equation with initial condition~\eqref{e:effe} and convolution kernel~\eqref{e:anisotropic}, when $\ee = 0.25$, $h = 0.01$. \label{fig:Test4}} \label{f:F}
\end{figure}
We now investigate the behavior of the support of the numerical solution in the case of anisotropic convolution kernels. More precisely, in the previous paragraph we have pointed out that a main drawback of the Lax-Friedrichs type method is that it does not preserve an important qualitative property of the analytic solutions of the nonlocal equations, namely the fact that, in the example we test at the previous paragraph, they satisfy~\eqref{e:rhoci}. Conversely, the Godunov type method preserves~\eqref{e:rhoci}. 

We now test this behavior in the case of a smoother initial datum than~\eqref{eq:init_Test3}. More precisely, we set 
\begin{equation} \label{e:effe}
    \bar \rho^{(F)}(x) = (x+1) \mathbbm{1}_{[-1,0]}(x)
\end{equation}
and we take the same convolution kernels as in~\eqref{e:anisotropic}. This implies that the solutions of the nonlocal Cauchy problem obtained by coupling~\eqref{e:effe} with~\eqref{e:nlburgers} satisfy~\eqref{e:rhoci}. In Test 5 we compute the numerical solutions with the Lax-Friedrichs type and the Godunov type methods and we display the snapshots of the solution in Figure~\ref{f:F}. Remarkably, as in Test 4 the solutions obtained by the Lax-Friedrichs type method do not satisfy~\eqref{e:rhoci}, while the solutions obtained by the Godunov type method satisfy~\eqref{e:rhoci}. 
Once again we find that the results obtained with the Lax-Friedrichs type method are not consistent with the analytic results and that the Godunov type scheme is more reliable.   
\subsection{Test 6 (Example C): positive density and isotropic convolution kernels}
	\label{ss:NumRes_Test2} 
	\begin{figure}[!h]
	\centering
	\includegraphics[width=0.8\textwidth]{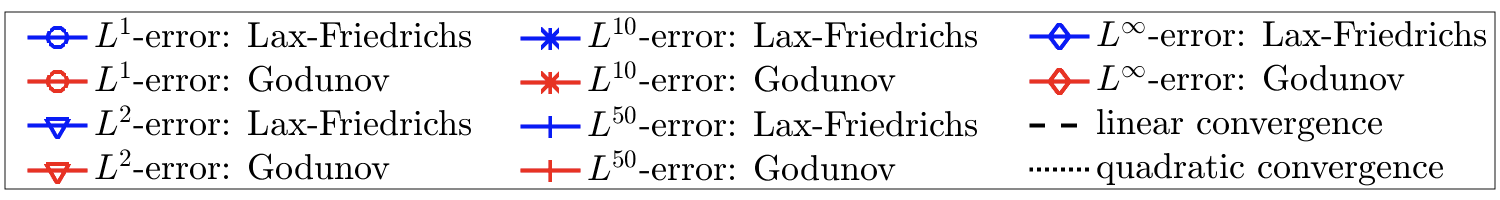}\\
	\includegraphics[width=0.48\textwidth]{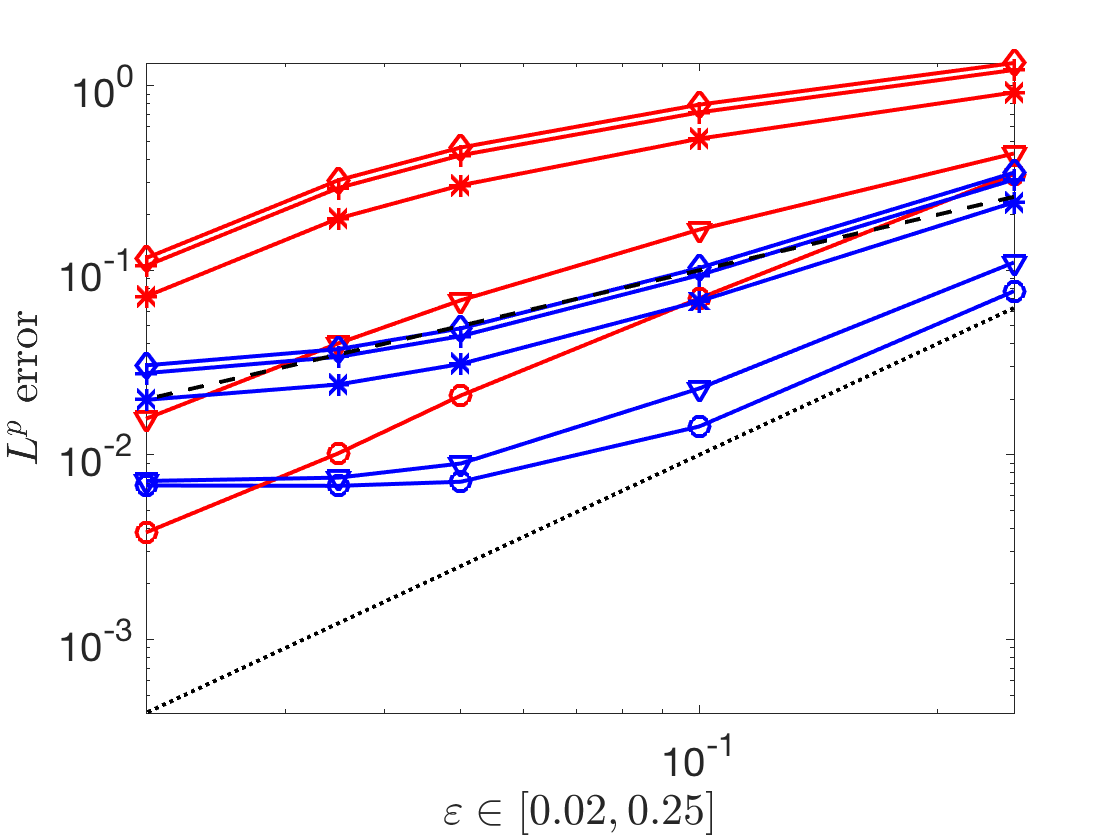} \hfill 
	\includegraphics[width=0.48\textwidth]{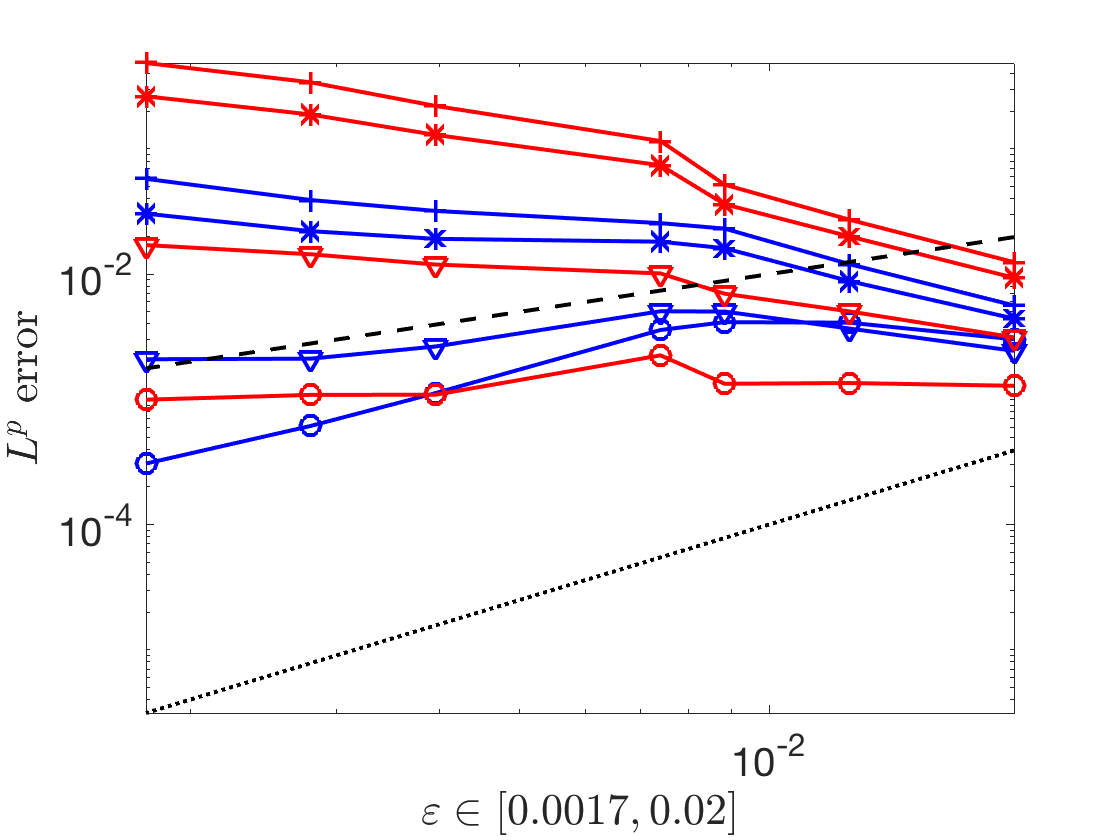}\\
	\hspace*{0.2\textwidth} (a) \hfill (b) \hspace*{0.2\textwidth} \textcolor{white}{.} 
	\caption{\footnotesize Test 6 (Example C), $L^p$-error at $t=2$, $p\geq 1$, at different values of $\ee$, comparing the nonlocal solution to the local solution for both numerical schemes: (a) fixed viscosity $h=0.001$, (b) varying viscosity $h=64 \ee^2$. \label{fig:Test2_conv-a}}
\end{figure}
\begin{figure}[!h]
	\centering 
	\footnotesize{Analytic solution} \\
	\includegraphics[width=\textwidth]{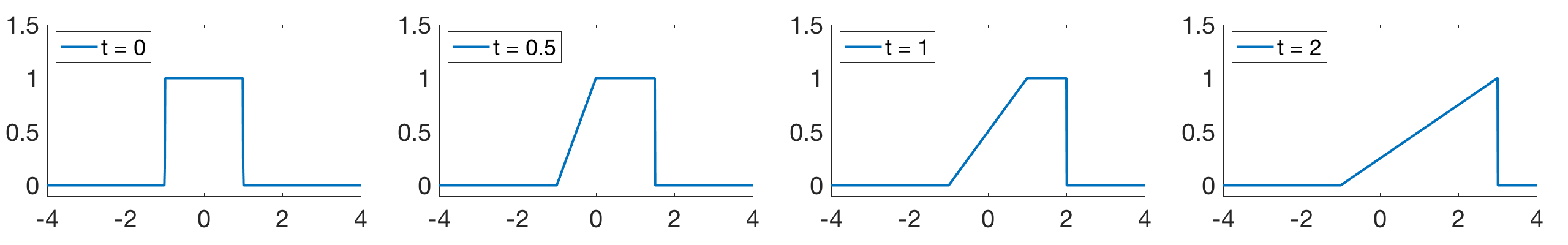} \\[1.em]
	\footnotesize{Numerical solution of the inviscid Burgers' equation with Lax-Friedrichs method} \\
	\includegraphics[width=\textwidth]{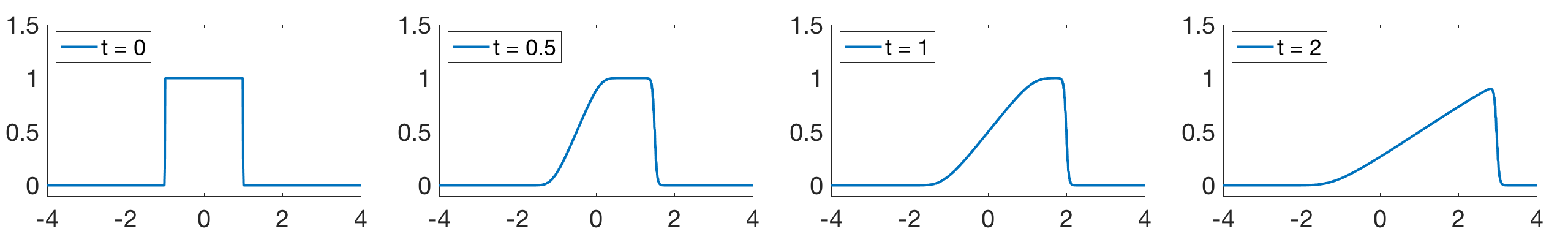} \\[1.em]
	\footnotesize{Numerical solution of the inviscid Burgers' equation with Godunov method} \\
	\includegraphics[width=\textwidth]{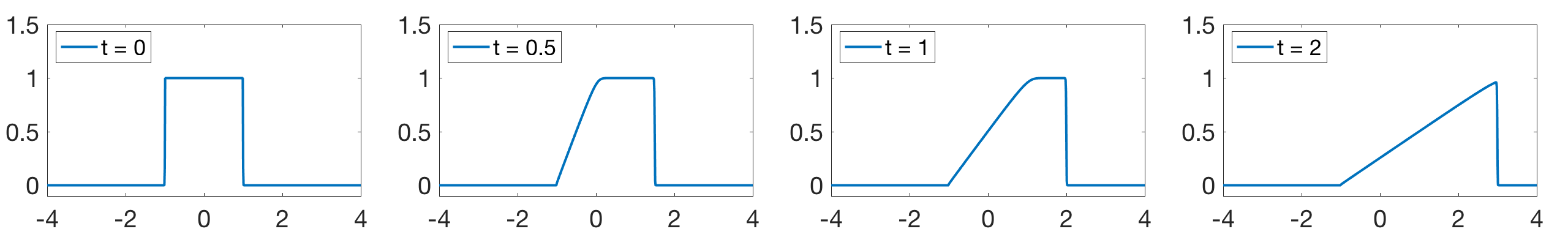} \\[1.em]
	\footnotesize{Numerical solution of the nonlocal Burgers' equation with Lax-Friedrichs method} \\
	\includegraphics[width=\textwidth]{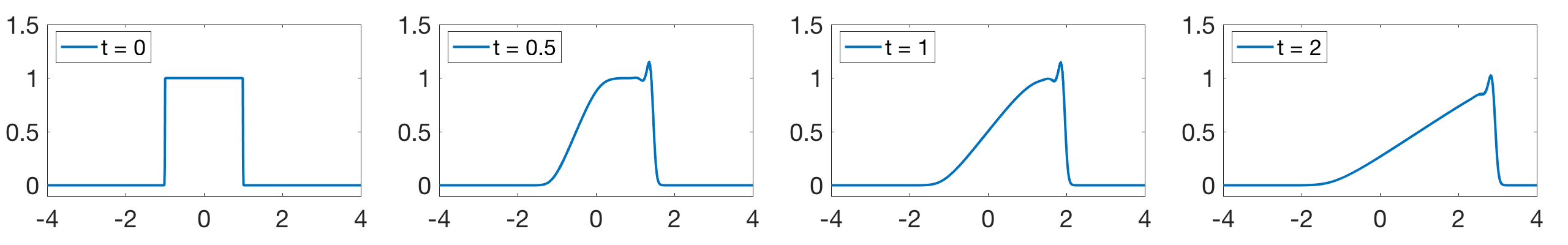} \\[1.em]
	\footnotesize{Numerical solution of the nonlocal Burgers' equation with Godunov method} \\
	\includegraphics[width=\textwidth]{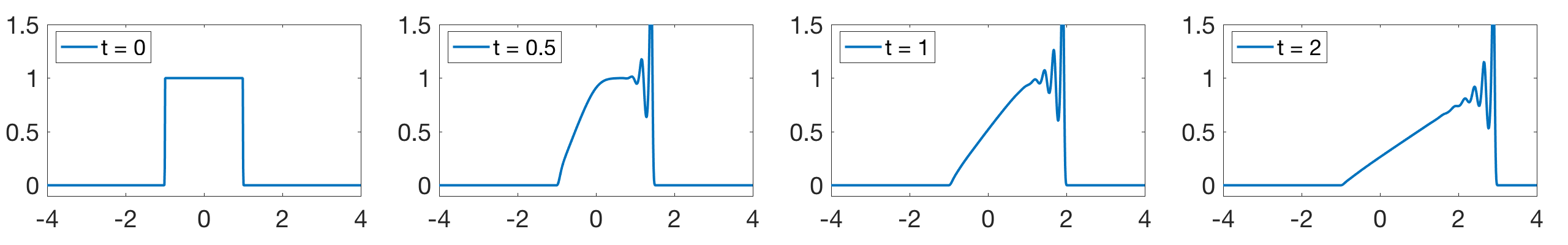} 
	\caption{\footnotesize Test 6 (Example C), snapshots of solution of Burgers' equation with initial condition and convolution kernel~\eqref{eq:init_Test2}, when $\ee = 0.2$, $h = 0.01$. \label{fig:Test2}}
\end{figure}
In Test 6 we take the same initial datum as in~\eqref{eq:init_Test2} and the same convolution kernels as in~\eqref{e:even}. Since the convolution kernels are even functions, we can apply the discussion in \S\ref{exC} and conclude that, for every $p>1$, the solutions of the nonlocal equations \emph{do not} converge to the entropy admissible solution of the (local) Burgers' equation strongly in $L^p$. 
 
In Test 6 we compute the numerical solution of the nonlocal equations with the Lax-Friedrichs type and with the Godunov type methods. We display the corresponding results in Figure~\ref{fig:Test2_conv-a} (convergence analysis) and in Figure~\ref{fig:Test2} (snapshots of the solution).  More precisely, Figure~\ref{fig:Test2_conv-a} shows the $L^p$ norm of the difference, computed at time $t=2$ between the numerical solutions of the nonlocal equations and the exact entropy solution of the (local) Burgers' equation. Here are the main comments. 
\begin{itemize}
\item[i)] In the results displayed in Figure~\ref{fig:Test2_conv-a}, part (a), we keep the space mesh $h$ fixed and we consider smaller and smaller values 
of the convolution parameter $\ee$. Both the Lax-Friedrichs and the Godunov type schemes results suggest that the $L^p$ norm converges to $0$, which is contradicted by the analytic results in~\S\ref{exC}. This is due to the fact that both schemes contain some numerical viscosity. 
\item[ii)] Figure~\ref{fig:Test2_conv-a}, part (b) displays numerical results where the space mesh $h$ is of the order of $\ee^2$. In this case the numerical results \emph{do not} suggest convergence in $L^p$ for $p>1$ and hence are consistent with the analytic results. This is most likely due to the fact that 
 the numerical viscosity  goes to $0$ very fast  and hence does not affect the investigation of the nonlocal-to-local limit. 
\end{itemize}
In a nutshell, Test 6 shows that the presence of the numerical viscosity jeopardizes the reliability of the numerical schemes and that the most reliable results are obtained by taking the smallest numerical viscosity. 
\subsection{Test 7 (Example E): isotropic convolution kernels and regular limit solution}
	\label{ss:Other_Test0}
	In Test 7 we consider the same monotone increasing initial datum $\bar{\rho}^{(E)}$ as in~\eqref{eq:init_TestE} and the same isotropic 
	convolution kernels as in~\eqref{e:even}. As pointed out in~\S\ref{sss:E}, owing to a result by Zumbrun~\cite[Proposition 4.1]{Zumbrun} we expect that in this case the solutions of the nonlocal equations converge to the (regular) solution of the Burgers' equation. 
	
\begin{figure}[!h]
	\centering
	\includegraphics[width=0.45\textwidth]{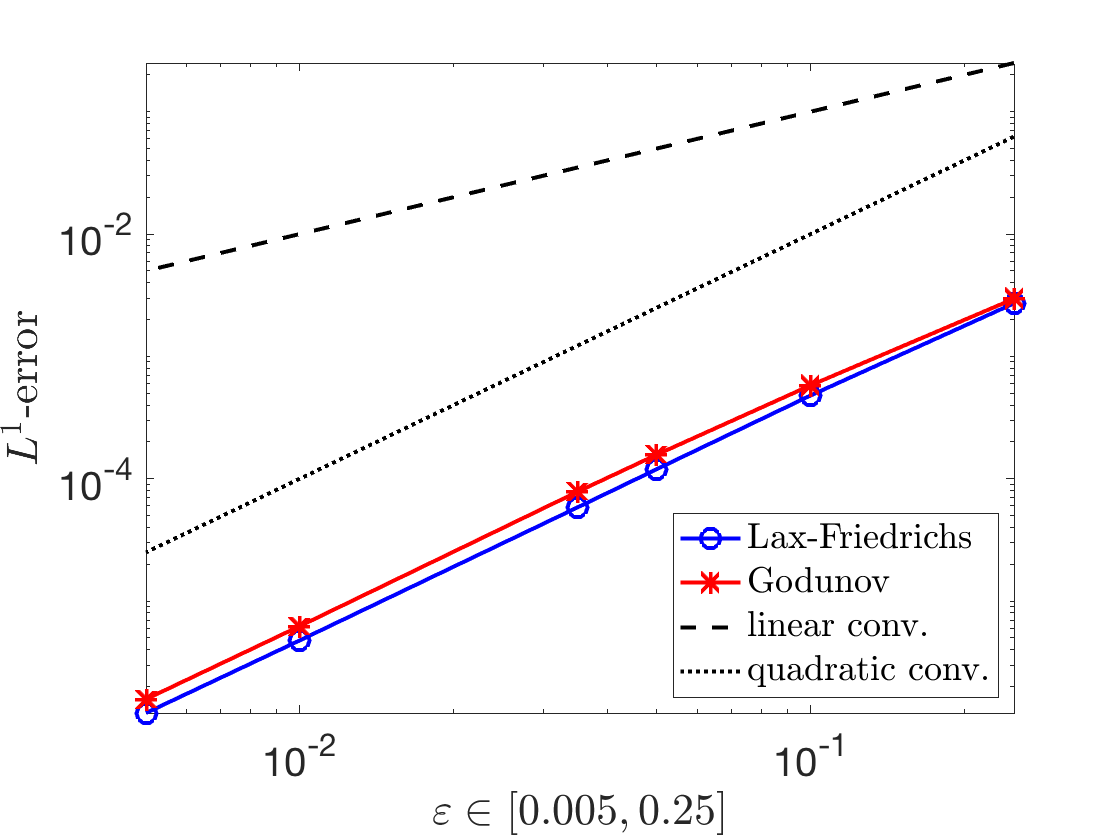} \hfill 
	\includegraphics[width=0.45\textwidth]{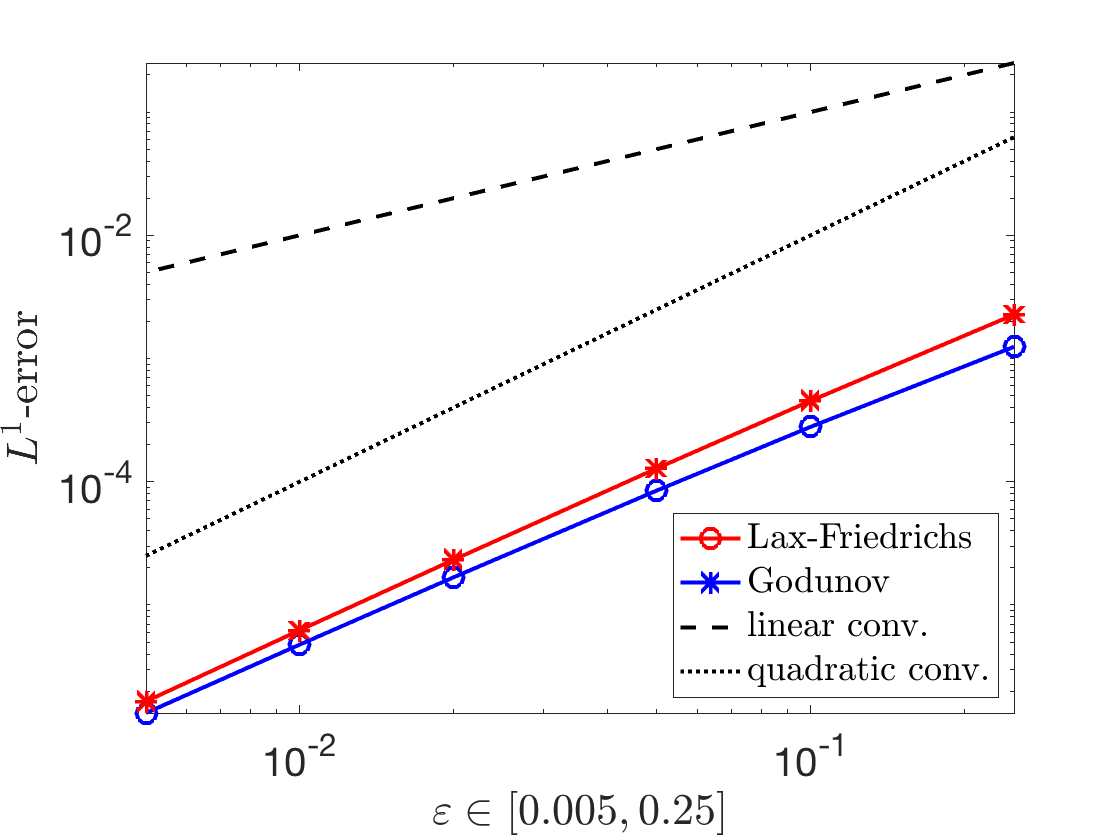}\\
	\hspace*{0.2\textwidth} (a) \hfill (b) \hspace*{0.2\textwidth} \textcolor{white}{.} 
	\caption{\footnotesize Test 7 (Example E), $L^1$-error at $t=2$, for different values of $\ee$, comparing the solutions of the nonlocal equations with the entropy solution of the (local) Burgers' equation for both numerical schemes, when (a) $h=0.001$, (b) $h=\ee /10$. \label{fig:Test0_conv-a}}
\end{figure}

	In Test 7 we compute the numerical solutions of the nonlocal equations with the Lax-Friedrichs type and the Godunov type methods. We display the results in Figure~\ref{fig:Test0_conv-a}: the results   
	obtained with both the Lax-Friedrichs and the Godunov type schemes suggest convergence. In the case of Test 7 this is consistent with the analytic results.  
\section{Conclusion}\label{s:end}
In this paper we have investigated the role of numerical viscosity in the study of the local limit of nonlocal conservation laws. We have shown that numerical viscosity is a very relevant feature as it severely affects the reliability of numerical methods. This claim is supported by the following instances: 
\begin{itemize}
\item Lax-Friedrichs type schemes have a very high numerical viscosity and erroneously suggest convergence in cases where convergence 
is ruled out by analytic considerations (see Tests 3, 4 and 6). Also, Lax-Friedrichs type scheme fail to capture relevant qualitative properties of the solutions of the nonlocal equations (see Tests 4 and 5). 
\item Godunov type schemes have lower numerical viscosity than Lax-Friedrichs type schemes and, at least in some cases, provide more reliable information on the nonlocal-to-local limit, see Tests 3 and 4. Also, they manage to capture relevant qualitative properties of the solutions of the nonlocal equations that are 
missed by the Lax-Friedrichs type schemes  (see Tests 4 and 5). 
\item The most reliable numerical results are obtained by using the Godunov type scheme and choosing a very low numerical viscosity (i.e. choosing the space mesh of the order $\ee^2$, where $\ee$ is the convolution parameter), see Tests 3 and 6. 
\end{itemize}
We feel that the present work could pave the way for several interesting developments. Note indeed that the Godunov type scheme provides more reliable results than the Lax-Friedrichs type scheme, but it is still not completely satisfactory and reliable for the analysis of the nonlocal-to-local limit. Also, it is more difficult to implement than the Lax-Friedrichs type scheme. It would be very interesting, for both numerical and analytic purposes
\footnote{Reliable numerical scheme could for instance provide valuable intuition on some analytic open questions. As an example, we mention traffic models with completely anisotropic convolution kernels, which take into account the fact that drivers only look forward, not backward, and hence decide their speed based on the downstream traffic density only. In this case Question~\ref{?} on the nonlocal-to-local limit is presently open, even if a recent counterexamples rules out the most ``natural'' strategy to establish convergence, see~\cite{CCS2}.   
 }, to introduce numerical schemes providing reliable results for the study of the local limit of nonlocal conservation laws.  
To this end, a possible future direction is working with Lax-Wendroff type schemes\footnote{We thank Giovanni Russo for this remark.}. This is motivated by the fact that the modified equation for the Lax-Wendroff equation is a third order equation with no viscous term, see~\cite{LeVeque}.
\section{Acknowledgments} The authors wish to thank Blanca Ayuso de Dios and Giovanni Russo for interesting discussions. 
GC is partially supported by the Swiss National Science Foundation grant 200021-140232 and by the ERC Starting Grant 676675 FLIRT. LVS is a member of the GNAMPA
group of INDAM and of the PRIN National Project ``Hyperbolic Systems of Conservation
Laws and Fluid Dynamics: Analysis and Applications''. 
MG was partially supported by the Swiss National Science Foundation grant P300P2-167681. 
Part of this work was done when MC and LVS were
visiting the University of Basel: its kind hospitality is gratefully acknowledged.
\bibliographystyle{plain}
\bibliography{singpar}

\end{document}